\documentclass[12pt,leqno]{article}
\usepackage[utf8]{inputenc}
\usepackage[official]{eurosym}
\usepackage{fullpage,theorem}
\usepackage{epic,eepic}
\usepackage{amssymb,amsmath}
\usepackage{accents}
\usepackage{cases}

\def\min{{\rm min}}
\def\max{{\rm max}}
\def\sup{{\rm sup}}
\def\otp{{\rm otp}}

\def\succ{{\rm Succ}}
\def\Lev{{\rm Lev}}
\def\dom{{\rm dom}}

\newtheorem{theorem}{Theorem}[section]
\newtheorem{lemma}[theorem]{Lemma}
\newtheorem{question}[theorem]{Question}
\newtheorem{corollary}[theorem]{Corollary}
\newtheorem{proposition}[theorem]{Proposition}
{\theorembodyfont{\rmfamily}
\newtheorem{definition}[theorem]{Definition}}
{\theorembodyfont{\rmfamily}
\newtheorem{example}[theorem]{Example}}
{\theorembodyfont{\rmfamily}
}
{\theorembodyfont{\rmfamily}

}
{\theorembodyfont{\rmfamily}
}
{\theorembodyfont{\rmfamily}
}

\title{
Prikry Forcing and Tree Prikry Forcing of Various Filters}
\author{Tom Benhamou\footnote{ Tel-Aviv University, Department of Mathematics, 
tombenhamou@mail.tau.ac.il}}
\date{\today}

\begin{document}
\maketitle
\begin{abstract}
    In this paper, we answer a question asked in \cite{MinimalPrikry} regarding a  Mathias criteria for Tree-Prikry forcing. Also we will investigate Prikry forcing using various filters. For completeness and self inclusion reasons, we will give proofs of many known theorems.    
\end{abstract}
\section{Introduction}

A well known result in abstract forcing theory, is that a forcing which preserves all cofinalities, also preserves cardinals. It is natural to consider the opposite question i.e. does every cardinal preserving forcing also preserves cofinalities?

 Karel Prikry introduced his classic Prikry forcing \cite{Prikry}, which was originally designated to give a counterexample for that statement i.e. a forcing notion, which preserves all cardinals and changes the cofinality of a measurable cardinal $\kappa$, while adding no new bounded subsets to $\kappa$.
The definition of Prikry forcing with a normal ultrafilter $U$ over $\kappa$, denoted by $P(U)$, uses the existence of a measurable cardinal.
Later it was shown
by Dodd and Jensen \cite{DoddJensen} that if there is a forcing notion which preserves all cardinals and changes cofinalities, there is an inner model with a measurable cardinal. 

The main feature of $P(U)$ is that a cofinal $\omega$-sequence is added to $\kappa$, while no new bounded subsets of $\kappa$ are added. Such sequences are usually called Prikry sequences. Mathias \cite{Mathias} found a criteria that ensures that an $\omega$-sequence is a Prikry sequence (see theorem \ref{Mathias}).

Devlin investigated a generalization of Prikry forcing over various kind of filters \cite{Devlin}. One of his results is the classification of filters for which this generalized version of Prikry forcing preserves cardinals and changes the cofinality of $\kappa$ to $\omega$. An important corollary of his work is that in order to preserve cardinals and change cofinalities, Prikry forcing with a non-normal filter can be used. Furthermore, Devlin and Paris proved \cite{Devlin} that the filters $U$, for which $P(U)$ adds no new bounded subsets to $\kappa$, are exactly Rowbottom ultrafilters (see definition \ref{VariousFilters}). In this paper we will determine which are the filters for which $P(U)$ has the Mathias criteria. 

\textbf{Theorem 3.21} $P(U)$ satisfies the Mathias criteria if and only if  $U$ is a Rowbottom ultrafilter.

Another interesting research in the field of Prikry forcing, is the investigation of intermediate ZFC models of Prikry forcing extensions.  Gitik, Koepke and Kanovei, proved that an intermediate ZFC model of Prikry forcing with a normal ultrafilter $U$, must also be a Prikry extension of the ground model for Prikry forcing with the same $U$ \cite{PrikryCase}.

In the absence of normality, there is a variation of Prikry forcing that can be defined--the Tree-Prikry forcing, which we denote by $P_T(\vec{U})$, where $\vec{U}$ is a tree of filters (see definition \ref{TreeOfFilters}). In Tree-Prikry extensions, as in regular Prikry extension, cardinals are preserved and an $\omega$-sequence is added to $\kappa$, while no new bounded subsets are added. Nonetheless, there are some differences. For example, the structure of intermediate models of generic extensions of $P_T(\vec{U})$ can be more complex. Gitik and the author proved that it is possible that a Cohen generic extension is an intermediate model of a Tree-Prikry forcing \cite{TomMoti}. This situation is not possible in the regular Prikry forcing (with a normal ultrafilter) according to the result of Gitik, Koepke and Kanovei. 

On the other hand, it is possible that a Tree-Prikry forcing extension is very simple. Koepke, Rasch and Schlicht proved that for certain trees of ultrafilters, the model obtained by the Tree-Prikry forcing is minimal i.e. has no proper intermediate extensions \cite{MinimalPrikry}. In their paper, they ask for a Mathias-like criteria for Tree-Prikry forcing. Such a criteria is given here: 

\textbf{Theorem 4.18} Let $\vec{U}$ be a tree of filters such that for every $a\in[\kappa]^{<\omega}$, $U_a$ is an ultrafilter which contains all the final segments. If $C\in[\kappa]^{\omega}$ is a sequence such that
\begin{enumerate}

\item $\sup(C)=\kappa$.
 \item For every $\langle A_a\mid a\in[\kappa]^{<\omega}\rangle$ such that $A_a\in U_a$, there exists $n<\omega$ such that for every $n\leq m<\omega$,
 $C(m)\in A_{C\restriction m}$.
 \end{enumerate}
 Then $C$ is Tree-Prikry-generic for $\vec{U}$.
\section{Some Theory of Filters}
The classical Prikry forcing can be performed using various types of filters, the filter combinatorical properties changes drastically the properties of the forcing, this discussion will take place in the next section. This section is devoted to basic definitions of filters and facts about them. 
\begin{definition}
A set $U\subseteq P(S)$ is a \textit{filter} over $S$ if:
\begin{enumerate}
    \item $\emptyset\notin U, \ S\in U$.
    \item $\forall A,B \ (A,B\in U\Rightarrow A\cap B\in U)$.
    \item $\forall A,B \ (A\in U\wedge A\subseteq B \Rightarrow B\in U)$.
    \end{enumerate}
Some additional definitions of types of filters:
\begin{enumerate}
\item $U$ is \textit{uniform} if for every $X\in U$, $|X|=|S|$.
\item $U$ is $\lambda$-\textit{complete} if for every $\beta<\lambda$, $\langle A_{\alpha}\mid \alpha<\beta\rangle$ such that $A_\alpha\in U$, $\cap_{\alpha<
\beta}A_\alpha\in U$.
\item $U$ is \textit{trivial} if there is $a\in S$ such that $\{a\}\in U$.
\item $U$ is an \textit{ultrafilter} if $\forall X\in P(S) ( X\in U\vee S\setminus X \in U)$.

\end{enumerate}
\end{definition}
Let us focus on filters over some regular cardinal $\kappa$. In the following definition and throughout the paper, $[A]^\alpha$ denotes the set of increasing sequences of elements of $A$ of order type $\alpha$, $[A]^{<\alpha}$ denotes increasing sequences of elements of $A$ of order type less than $\alpha$. Every set of ordinals is ordered naturally by the
usual order of ordinals. We identify a set of ordinals with the  sequence of its natural increasing enumeration. 

Let $f:A\rightarrow B$ be a function, denote the \textit{point-wise image} of a set $X\subseteq A$ by $$f[X]=\{f(x)\mid x\in X\}.$$  Also, define the \textit{preimage} of a set $Y\subseteq B$ by $$f^{-1}[Y]=\{x\in X\mid f(x)\in Y\}.$$

\begin{definition}\label{VariousFilters}
Let $U$ be a filter over $\kappa$.
\begin{enumerate}
\item $U$ is \textit{normal} if for any $A\in U$ and any $F:A\rightarrow\kappa$ such that for all $\alpha\in A$, $F(\alpha)<\alpha$, there is $A'\in P(A)\cap U$ such that $F\restriction A'$ is constant.
\item $U$ is \textit{Rowbottom} if for any $A\in U$ and any $F:[A]^{<\omega}\rightarrow X$ such that $|X|<\kappa$, there is $A'\in P(A)\cap U$ such that for every $n<\omega$, $F\restriction [A']^n$ is constant.
\item $U$ is \textit{shrinking} if for any $A\in U$ and any $F:[A]^{<\omega}\rightarrow \lambda^+$ such that $\lambda^+<\kappa$, there is $A'\in P(A)\cap U$ such that  $|F[[A']^{<\omega}]|\leq\lambda$.
\item $U$ is \textit{Ramsey} if for any $A\in U$ and any $F:[A]^2\rightarrow\{0,1\}$,
there is $A'\in P(A)\cap U$ such that $F\restriction[A']^2$ is constant.
\end{enumerate}
\end{definition}

A measurable cardinal is a cardinal which carries a $\kappa$-complete non-trivial ultrafilter. It is well known that the existence of a measurable cardinal is not provable from ZFC.
\begin{proposition}\label{BasicProp}
Let $\kappa$ be weakly inaccessible cardinal. For every non-trivial uniform filter $U$ over $\kappa$:
$$U \ is \ normal \Rightarrow \ U\ is \ Rowbottom \Rightarrow \ U \ is \ shrinking \ and \ Ramsey.$$
\end{proposition}
\textit{Proof}: For the implication ``normal $\Rightarrow$ Rowbottom", see \cite{Jech} theorem 10.22. The implication  ``Rowbottom $\Rightarrow$ shrinking and Ramsey" follows directly from the definitions. $\square$

It is known that the implications in proposition \ref{BasicProp} cannot be reversed \cite{Devlin2}.
\begin{definition}
Let $U$ be an ultrafilter over $S$ and $f:S\rightarrow X$. The \textit{Rudin-Keisler projection of $U$ by $f$} is an ultrafilter in $X$ defined by $f_*(U)=\{Y\subseteq X\mid f^{-1} [Y]\in U\}$. For ultrafilters $U$ and $W$, define: 
\begin{enumerate}
    \item $W\leq_{R-K} U$ if $\exists f \ W=f_*(U)$.
    \item $U\equiv_{R-K} W$ if $U\leq_{R-K} W \wedge W\leq_{R-K} U$.
\end{enumerate}
\end{definition}
\begin{proposition}\label{RKEquivalence}
 Let $U$ be an ultrafilter over $I$ and $W$ an ultrafilter over $J$. Assume that $|I|=|J|$, then $U\equiv_{R-K} W$ iff there is a bijection $f:J\rightarrow I$ such that $f_*(W)=U$.
\end{proposition}
\textit{Proof}: Assume that $f$ is a bijection such that $f_*(W)=U$, then $U\leq_{R-K}W$. $f$ is invertible and $f^{-1}_*(U)=W$ witnessing $W\leq_{R-K} U$. 
For the other direction, assume that $U\equiv_{R-K}W$, then there are $f:J\rightarrow I$ and $g:I\rightarrow J$ such that $f_*(W)=U$ and $g_*(U)=W$. It follows that $$(f\circ g)_*(U)=f_*(g_*(U))=U.$$ 
Define $h=f\circ g$, then $h:I\rightarrow I$ and for every $X\subseteq I$, $X\in U\Leftrightarrow h[X]\in U$. We claim that $$\{x\in I\mid h(x)=x\}\in U.$$ Otherwise, since $U$ is an ultrafilter $B=\{x\in I\mid h(x)\neq x\}\in U$.
Assume that $\{ b_i\mid i< |B|\}$ is some enumeration of $B$, define $\phi: B\rightarrow \{0,1,2\}$ by recursion, such that $$\forall i<|B| \ \phi(b_i)\neq\phi(h(b_i)).$$ 
Let $\phi(b_0)=0$, assume that $\phi\restriction\{b_i\mid i<j\}$ is defined for some $j<|B|$. Also assume that for every $k<\omega$ and $i<j$, either $h^k(b_i)\notin B$ or $\phi(h^k(b_i))$ is defined ($h^k$ denoted the composition of $h$, $k$ times). 
If there are $m,k<\omega$ and $i<j$ such that $h^m(b_j)=h^k(b_i)$, take minimal such $m$. It follows that for every $m\leq n<\omega$, either $h^n(b_j)\notin B$ or $\phi(h^n(b_j))$ is defined. 
Assume for example that $\phi(h^m(b_j))=2$, define  $$\phi(h^k(b_j))=k \ mod(2), \text{ for }k<m.$$ If there are no $m,k,i$ such that $h^m(b_j)=h^k(b_i)$, again we separate into two cases.
First, if $$\exists k_1<k_2<\omega \ h^{k_1}(b_j)=h^{k_2}(b_j),$$ find minimal such $k_2$ and define for $m<k_2-1$, $\phi(h^m(b_j))=m \ mod(2)$ and  $\phi(h^{k_2-1}(b_j))=2$. 
Finally, if there are no $k_1,k_2<\omega$ such that $h^{k_1}(b_j)=h^{k_2}(b_j)$, for every $k$ such that $h^k(b_j)\in B$ define $\phi(h^k(b_j))=k\ mod 2$. 
For $i\in\{0,1,2\}$, consider $B_i=\phi^{-1}[\{i\}]\subseteq B$, then $$B_1\uplus B_2\uplus B_3=B\in U.$$ There is exactly one $i\in\{1,2,3\}$ such that $B_i\in U$.
Without loss of generality, suppose that $B_1\in U$ and $B_2\cup B_3\notin U$. Then $h[B_1]\in U$ and $h[B_1]\cap B\in U$. For every $h(b)\in h[B_1]\cap B$, it follows by definition that $\phi(h(b))\neq \phi(b)=1$, implying that $h(b)\in B_2\cup B_3$. 
Hence $$h[B_1]\cap B\subseteq B_2\cup B_3,$$ and so $B_2\cup B_3\in U$, contradiction. Therefore $A=\{x\in I\mid h(x)=x\}\in U$ and $g$ is 1-1 on this $U$-large set. 
Decompose $I,J$ each into two sets $I_1,I_2$ and $J_1,J_2$ respectively, such that $$|I_1|=|I_2|=|J_1|=|J_2|=|I|=|J|.$$ Exactly one of $I_1,I_2$ is in $U$ and one of $J_1, J_2$ is in $W$. Suppose for example that $I_1\in U$ and $J_1\in W$. Since $g_*(U)=W$, $g^{-1}[J_1]\in U$. Define $\rho:I\rightarrow J$, first $$\rho\restriction (I_1\cap A\cap  g^{-1}[J_1])=g\restriction (I_1\cap A\cap  g^{-1}[J_1]).$$ Note that $$I_2\subseteq I\setminus (I_1\cap A\cap  g^{-1}[J_1]) \text{ and } J_2\subseteq J\setminus g[I_1\cap A\cap  g^{-1}[J_1]].$$
So $$|I\setminus (I_1\cap A\cap  g^{-1}[J_1])|=|I|=|J|=|J\setminus g[ I_1\cap A\cap  g^{-1}[J_1]]|.$$ Therefore, let $\rho\restriction (I\setminus (I_1\cap A\cap  g^{-1}[J_1]))$ be any bijection with $J\setminus g[ I_1\cap A\cap  g^{-1}[J_1]]$. Then $\rho:I\rightarrow J$ is a bijection. Moreover, $\rho_*(U)=W$ since
$$X\in \rho_*(U)\Leftrightarrow \rho^{-1}(X)\in U\Leftrightarrow \rho^{-1}(X)\cap I_1\cap A\cap  g^{-1}[J_1]\in U\Leftrightarrow$$
$$g[\rho^{-1}[X]\cap I_1\cap A\cap  g^{-1}[J_1]]\in W\Leftrightarrow X\cap g[I_1\cap A]\cap J_1\in W\Leftrightarrow X\in W.$$
$\square$
\begin{proposition}\label{Rownor}
 The following are equivalent for every non-trivial $\kappa$-complete
ultrafilter $U$ over $\kappa$:
\begin{enumerate}
\item $U$ is $\leq_{R-K}$-minimal among $\kappa$-complete ultrafilters.
\item $U\equiv_{R-K}W$ for some normal ultrafilter $W$.
\item $U$ is Rowbottom.
\item $U$ is Ramsey.
\end{enumerate}
\end{proposition}
\textit{Proof}: For
 $(1)\rightarrow(2)$, let $Ult(V,U)$ be the ultrapower of $V$ by the ultrafilter $U$. Let $M$ be the transitive collapse of $Ult(V,U)$. Consider the corresponding elementary embedding  $$j_U:V\rightarrow Ult(V,U)\simeq M,$$ with critical point $\kappa$. 
 Let $\pi:\kappa\rightarrow\kappa$ be a function representing $\kappa$ in $Ult(V,U)$ i.e. $[\pi]_U=\kappa$. Then $\pi_*(U)\leq_{R-K} U$ and by the minimality assumption of $U$, $U\equiv_{R-K} \pi_*(U)$. 
 To see that $\pi_*(U)$ is normal, let $A\in \pi_*(U)$ and $f:A\rightarrow\kappa$, such that $\forall\alpha\in A$, $f(\alpha)<\alpha$. By definition of $\pi_*(U)$, $\pi^{-1}[A]\in U$ and for every $\beta\in \pi^{-1}[A]$, $f(\pi(\beta))<\pi(\beta)$. It follows that $$\{\beta<\kappa\mid f(\pi(\beta))<\pi(\beta)\}\in U.$$ 
 In terms of the ultrapower, $[f\circ\pi]_U<[\pi]_U=\kappa$. In particular, there is $\gamma<\kappa$, $$[f\circ\pi]_U=\gamma=j_U(\gamma)=[c_\gamma]_U,$$ where $c_\gamma$ is the constant function with value $\gamma$. Therefore, $D=\{\alpha<\kappa\mid f(\pi(\alpha))=\gamma\}\in U$, $\pi[D]\in \pi_*(U)$ and $f\restriction \pi[D]$ is constantly $\gamma$. 
 For $(2)\rightarrow(3)$, assume that $g:\kappa\rightarrow\kappa$ witnesses $U\equiv_{R-K}W$, for some normal ultrafilter $W$. Let $A\in U$ and $F:[A]^{<\omega}\rightarrow X$, where $|X|<\kappa$. 
 Define $G:[g^{-1}[A]]^{<\omega}\rightarrow[A]^{<\omega}$ by $$G(\langle\alpha_1,...,\alpha_n\rangle)=\langle g(\alpha_1),...,g(\alpha_n)\rangle.$$
 Define $F^*:[g^{-1}[A]]^{<\omega}\rightarrow X$ by $F^*=F\circ G$. Since $U=g_*(W)$,  $g^{-1}[A]\in W$ and by normality of $W$, there is $B\in W$, such that for every $n<\omega$, $F^*\restriction [B]^n$ is constant. 
 Let $A'=g[B]$, then $A'\in U$. For every $n<\omega$ and $\vec{\alpha},\vec{\beta}\in [A']^n$, there are $\vec{\alpha}',\vec{\beta}'\in [B]^n$ such that $G(\vec{\alpha}')=\vec{\alpha}$ and $G(\vec{\beta}')=\vec{\beta}$. It follows that
 $$F(\vec{\alpha})=F(G(\vec{\alpha}'))=F^*(\vec{\alpha}')=F^*(\vec{\beta}')=F(G(\vec{\beta}'))=F(\vec{\beta}).$$
 Hence $F\restriction[A']^n$ is constant. For
$(3)\rightarrow(4)$, use proposition \ref{BasicProp}. Finally
$(4)\rightarrow(1)$, let $V$ be a $\kappa$-complete ultrafilter such that  $V\leq_{R-K} U$ and $ \neg (V\equiv_{R-K} U)$. To see that $V$ is trivial,
consider $f:\kappa\rightarrow\kappa$ such that $f_*(U)=V$. Define $g:[\kappa]^2\rightarrow\{0,1\}$ by $$g(\langle\alpha,\beta\rangle)=1\Leftrightarrow f(\alpha)=f(\beta).$$ There is $A\in U$ such that $g\restriction[A]^2$ is constant. 
If this constant is $1$, then $f$ is constant on $A$, say with value $\gamma$. In particular $\{\gamma\}=f[A]\in V$, implying $V$ is trivial. If this constant is $0$ then $f$ is one to one on $A$. As in proposition \ref{RKEquivalence}, this implies $V\equiv_{R-K}U$, contradiction.
$\square$

For more information about filters see \cite{Kanamori}, \cite{Ketonen} or \cite{Devlin2}.
\section{Basic Prikry Forcing with Filters}

In this paper, we assume that the reader is familiar with forcing theory. We follow similar notations to \cite{Cumm} and use the Jerusalem-style notation of order i.e. $p\leq q$ means that ``$q$ is stronger than $p$". Accordingly, let $\mathbb{P}$ be a forcing notion, then the weakest condition in $\mathbb{P}$ is denoted by $0_{\mathbb{P}}$.  The ground model will be denoted by $V$ and if $\mathbb{P}\in V$ is a forcing notion, then $p\Vdash_{\mathbb{P}}\sigma$ is the statement ``$p$ forces $\sigma$". $\Vdash_{\mathbb{P}}\sigma$ means $0_{\mathbb{P}}\Vdash\sigma$. $\mathbb{P}$-names will be denoted with a dot accent e.g. $\dot{x}, \ \dot{f}, \ \dot{\gamma}$, etc. 
If $x\in V$, we will abuse notation by allowing $x$ to appear in a statements of the forcing language, where formally the canonical name of $x$ should have appeared.  If $\dot{x}$ is a name such that $\Vdash_{\mathbb{P}}\dot{x}\in V$, then the statement $p || \dot{x}$ means ``$\exists y\in V$ such that $p\Vdash_{\mathbb{P}}\dot{x}=y$". 
 For general information about forcing we refer the reader to \cite{Kunen} or \cite{Jech}.
 
 Let us introduce the Prikry forcing with a filter $U$ over a weakly inaccessible cardinal. 
\begin{definition} Let $U$ be a filter over a weakly inaccessible $\kappa$. The underlining set of the \textit{Prikry forcing}, denoted by $P(U)$, is the set of all elements of the form $\langle t_1,...,t_n, A\rangle$, where $\langle t_1,...,t_n\rangle\in[\kappa]^{<\omega}$, $A\in U$ and $\min(A)>t_n$. For $p=\langle t_1,...,t_n,A\rangle$, $q=\langle s_1,...,s_m,B\rangle$, define $p\leq q$ iff:
\begin{enumerate}
\item $n\leq m$.
\item $\forall i\leq n$, $t_i=s_i$.
\item $s_{n+1},...,s_m\in A$.
\item $B\subseteq A$.
\end{enumerate}
There is an important suborder $\leq^*\subseteq \leq$, defined by $p\leq^*q$ iff $p\leq q \wedge n=m$. \end{definition}

Some additional notations will be used. Let $p=\langle t_1,...,t_n,A\rangle\in P(U)$, define:
 \begin{enumerate}
 \item $t(p)=\langle t_1,...,t_n\rangle$, $n(p)=n$, $A(p)=A$.
 \item For $i\leq n$, $t_i(p)=t_i$.
 \item If $t,s\in[\kappa]^{<\omega}$,
 $s^{\frown}t$ denotes the concatenation of these sequences.
 \end{enumerate}
 
 For $\vec{\alpha}=\langle\alpha_1,...,\alpha_m\rangle\in [A]^{<\omega}$, $A'\subseteq A(p)$ and $A'\in U$, define:
 \begin{enumerate}
 \item   $p^{\frown}\langle \vec{\alpha},A'\rangle= \langle t_1,...,t_m,\alpha_1,...,\alpha_m,A'\rangle$.
\item $p^{\frown}\vec{\alpha}=p^{\frown}\langle\vec{\alpha},A(p)\setminus (\max(\vec{\alpha})+1)\rangle$.
 \end{enumerate}
\begin{lemma}\label{PrikryProp} For any filter $U$ over $\kappa$:
\begin{enumerate}
\item $P(U)$ is $\kappa^+$-Knaster. 
\item If $U$ is $\lambda$-complete, then $\leq^*$ is $\lambda$-closed.
\end{enumerate}
\end{lemma}
\textit{Proof}: For $(1)$, given any $\langle p_i\mid i<\kappa^+\rangle$, there is $A\subseteq \kappa^+$ of cardinality $\kappa^+$, $n^*<\omega$ and $t_1,...,t_{n^*}$  such that $$\forall i\in A , \ n(p_i)=n^*, \ t(p_i)=\langle t_{1},...,t_{n^*}\rangle.$$ Let $i,j\in A$, since $U$ is a filter, $A(p_i)\cap A(p_j)\in U$. It follows that $$\langle t_1,..,t_{n^*}, A(p_i)\cap A(p_j)\rangle\in P(U)$$ is a common extension of $p_i$ and $p_j$.
For $(2)$, let $\langle p_i\mid i<\beta\rangle$ where $\beta<\lambda$, be an $\leq^*$-increasing sequence. Since $U$ is $\lambda$-complete, $$A^*:=\cap_{i<\beta}A(p_i)\in U , \  \langle  t(p_0),A^*\rangle\in P(U).$$ The condition $\langle  t(p_0),A^*\rangle$ is an upper bound to the sequence $\langle p_i\mid i<\beta\rangle$.
$\square$
\begin{example}\label{pathologies}
The following simple examples suggests that without further assumptions about the filter, $P(U)$ might be degenerate. 
\begin{enumerate}
\item If $U$ is trivial, then $P(U)$ is atomic i.e. for all $p\in P(U)$ there is an atom $p\leq a$.\\
\textit{Proof}:
Suppose $U$ is trivial, then there is $\alpha<\kappa$ such that $\{\alpha\}\in U$. It follows that for every $p\in P(U)$, there is an atom $p\leq^*p^{\frown}\langle\{\alpha\}\rangle$.
\item Assume $\lambda$ is measurable, $U'$ is a normal ultrafilter over $\lambda$. Let $\lambda<\kappa$ be some regular cardinal, define $U=\{X\subseteq\kappa\mid X\cap\lambda\in U'\}$, then $\Vdash_{P(U)} \kappa$ is regular.\\
\textit{Proof}: Let $p\in P(U)$, define $$p_\lambda:= p^{\frown}\langle A(p)\cap\lambda\rangle\in P(U').$$ Note that $p_\lambda\geq p$, hence $P(U')$ is dense subset of $P(U)$. It follows that forcing with $P(U)$ is the same as forcing with $P(U')$. By lemma \ref{BasicProp}, $P(U')$ is $\lambda^+$-Knaster and $\kappa\geq\lambda^+$, so $\Vdash_{P(U')} \kappa$ is regular. 
\end{enumerate}
\end{example}
In order for $P(U)$ to be in the spirit of Prikry, it should at least change the cofinality of $\kappa$ to $\omega$. By the previous examples, we should at least assume that $U$ contains all final segments  $(\lambda,\kappa):=\kappa\setminus(\lambda+1)$. Since $\kappa$ is regular, such a filter must also be non-trivial and uniform. 

The next definition is of a forcing which is forcing equivalent to $P(U)$ and will be easier to manage.
\begin{definition}
The underlining set of the forcing $P^*(U)$ is the set of all $f:X\rightarrow\{0,1\}$ such that $X\subseteq\kappa$, $|f^{-1}[\{1\}]|<\omega$ and $\kappa\setminus X\in U$. The order is defined by $f\leq g$ iff $f\subseteq g$.
\end{definition}
\begin{lemma}\label{Iso}
Let $U$ be a filter over $\kappa$ that contains all final segments, then:
\begin{enumerate}
\item $P^*(U)$ is forcing equivalent to $P(U)$.
\item If $U\equiv_{R-K} W$ then $P^*(U)\simeq P^*(W)$ in a canonical way.
\end{enumerate}
\end{lemma}
\textit{Proof}:
For $(1)$, we will construct a dense complete embedding $i:P(U)\rightarrow P^*(U)$
(For the definition of complete embedding see for example \cite{Kunen} or \cite{Shelah}). Define $i(p)=f$, where
$$\dom(f)=\kappa\setminus A(p), \ f(\alpha)=1\Leftrightarrow\alpha\in t(p).$$
Since $p\in P(U)$, $i(p)\in P^*(U)$. 
If $p\leq q$ then $$\dom(i(p))=\kappa\setminus A(p)\subseteq \kappa\setminus A(q),$$ and for $\alpha\in \dom(i(p))$,
$$i(p)(\alpha)=1\Leftrightarrow \alpha\in t(p)\Leftrightarrow\alpha\in t(q)\Leftrightarrow i(q)(\alpha)=1,$$ 
hence $i(p)\leq i(q)$. Assume that $i(p_1),i(p_2)$ are compatible, then $f=i(p_1)\cup i(p_2)$ most be a function and $i(p_1),i(p_2)\leq f$.
Without loss of generality, assume that $$(i(p_1))^{-1}[\{1\}]\cap\max(t(p_2))+1=(i(p_2))^{-1}[\{1\}].$$ 
By definition, $t(p_k)=(i(p_k))^{-1}[\{1\}]$ for $k\in\{1,2\}$. Hence $t(p_1)\cap\max(t(p_2))+1=t(p_2)$.
Moreover, $$t(p_1)\setminus t(p_2)=i(p_1)^{-1}[\{1\}]\setminus i(p_2)^{-1}(\{1\})\subseteq \kappa\setminus(\kappa\setminus A(p_1))=A(p_1).$$ 
Thus $p_1^{\frown}\langle A(p_1)\cap A(p_2)\rangle\in P(U)$ is a common extension of $p_1$ and $p_2$.
To see that $i[P(U)]$ is dense in $P^*(U)$, let $f\in P^*(U)$ and enumerate $f^{-1}[\{1\}]=\langle\alpha_1,...,\alpha_n\rangle$.
Since $U$ contains all final segments,  $A^*:=(\kappa\setminus \dom(f))\cap(\alpha_n,\kappa)\in U$.
Define $$p=\langle \alpha_1,...,\alpha_n,A^*\rangle\in P(U),$$
then $\dom(i(p))=\dom(f)\cup\alpha_n\supseteq \dom(f)$ and if $f(\alpha)=1$, then  $\alpha=\alpha_j$ for some $1\leq j\leq n$. By definition, $i(p)(\alpha)=1$, thus $i(p)\geq f$.
For $(2)$, let $f:\kappa\rightarrow\kappa$ be a bijection witnessing $U\equiv_{R-K} W$ i.e. $f_*(U)=W$. Construct $F:P^*(U)\rightarrow P^*(W)$,
$$\dom(F(g))=f[\dom(g)], \  F(g)(x)=g(f^{-1}(x)).$$ Note that $\kappa\setminus \dom(F(g))=f[\kappa\setminus \dom(g)] \in W$. Also $F(g)^{-1}[\{1\}]=f[g^{-1}[\{1\}]]$, therefore $|F(g)^{-1}[\{1\}]|<\omega$. Clearly, since $f$ is invertible, $F$ is  invertible. The proof that $F$ is order preserving, is a straightforward verification.  
$\square$
\begin{definition}
Let $G$ be $P(U)$-generic. 
Define the \textit{Prikry-generic sequence} derived from $G$ by $$C_G=\bigcup\{t(p)\mid p\in G\}.$$ A sequence $C$ is \textit{Prikry-generic for} $U$, if $C=C_G$ for some generic $G\subseteq P(U)$. 
\end{definition}

\begin{lemma}\label{Genericseq}
 Let $U$ be a filter which contains all the final segments, let $G$ be generic for $P(U)$, then:
\begin{enumerate}
\item $\sup(C_G)=\kappa$.
\item $\otp(C_G)=\omega$.
\item $\forall A\in U, \ C_G\setminus A$ is finite.
\item Suppose that $f:\kappa\rightarrow \kappa$ witnesses $U\equiv_{R-K} W$. If $C_G$ is Prikry-generic for $U$, then $f[C_G]$ is Prikry-generic for $W$.
\end{enumerate}
\end{lemma}
\textit{Proof}: For $(1)$, let $\beta<\kappa$ and $p\in P(U)$. Pick $\beta<\beta_1$ such that $A(p)\cap(\beta,\beta_1)\neq\emptyset$, 
and let $r\in A(p)$ be such that $\beta<r<\beta_1$, then  $$(\beta_1,\kappa)\in U\Rightarrow A(p)\cap (\beta_1, \kappa)\in U\Rightarrow p^{\frown}\langle r, A(p)\cap (\beta_1,\kappa)\rangle\in P(U).$$
Therefore, $p\leq p^{\frown}\langle r, A(p)\cap (\beta_1,\kappa)\rangle$ and $$p^{\frown}\langle r, A(p)\cap (\beta_1,\kappa)\rangle\Vdash \sup(\dot{C_G})>\beta.$$
By density, there is a condition $q_\beta\in G$ such that $q_\beta\Vdash \sup(\dot{C_G})>\beta$. It follows that $\sup(C_G)>\beta$, this is true for every $\beta<\kappa$, hence $\sup(C_G)=\kappa$.
For $(2)$, fix $\delta<\kappa$, for every $p\in P(U)$ define $p_\delta:=p^{\frown}\langle[\delta,\kappa)\cap A(p)\rangle\in P(U)$. Then $p\leq^* p_\delta$ and $$p_\delta\Vdash |\dot{C_G}\cap \delta|=n(p)<\omega.$$ By density, for every $\delta<\kappa$ there is such $p_\delta\in G$, hence $|C_G\cap\delta|=n(p_\delta)$. It follows that $\otp(C_G)\leq\omega$, otherwise, there is $\delta<\kappa$ such that $\otp(C_G\cap\delta)=\omega$, contradiction. If $|C_G|<\omega$, then $\alpha=\max(C_G)<\kappa$, this contradicts $(1)$. We conclude that $\otp(C_G)=\omega$.
For $(3)$, let $A\in U$ and $p\in P(U)$. There is $p^*=p^{\frown}\langle A\cap A(p)\rangle\in P(U)$, such that $p\leq^* p^*$.
We claim that $$p^*\Vdash \dot{C_G}\setminus \max(t(p^*))+1\subseteq A.$$ To see this,
let $H$ be generic with $p^*\in H$ and let $\alpha\in C_H\setminus (\max(t(p^*))+1)$.
There is $r\in H$ such that $r\geq p^*$ and $\alpha\in t(r)$.
By definition, $$\alpha\in t(r)\setminus t(p^*)\subseteq A(p^*)\subseteq A,$$ thus $C_H\setminus\max(t(p^*))+1\subseteq A$.
By density, there is $p^*\in G$ such that  $$p^*\Vdash\dot{C_G}\setminus \max(t(p^*))+1\subseteq A.$$
Therefore, $C_G\setminus A\subseteq \max(t(p^*))$, by $(1),(2)$ it follows that $C_G\setminus A$ is finite.
Finally for $(4)$, recall the definition of $F$ and $i$ from lemma \ref{Iso}.
Let $G$ be $P(U)$-generic, then $H=i^{-1}[F[i[G]]]$ is $P(W)$-generic.
Let $\alpha\in C_G$, find $p\in G$ with $\alpha\in t(p)$, then $$\alpha\in t(p)\Leftrightarrow i(p)(\alpha)=1\Leftrightarrow F(i(p))(f(\alpha))=1.$$ Since $i$ is a dense embedding, there is $q\in H$ such that $i(q)\geq F(i(p))$. Therefore, $f(\alpha)\in C_H$ and $f[C_G]\subseteq C_H$. The other inclusion is symmetric, concluding $f[C_G]=C_H$.
$\square$
\begin{corollary}
$cf^{V[G]}(\kappa)=\omega$.
\end{corollary}
\begin{definition}
Let $C\in[\kappa]^{\omega}$, define the \textit{filter generated by $C$}, \begin{center}
$G_C=\{p\in P(U)\mid t(p)=C\cap (\max(t(p))+1)\wedge C\setminus (\max(t(p))+1)\subseteq
A(p)\}$.
\end{center}
\end{definition}
\begin{lemma}\label{Filter}
For every $C$, $G_C$ is a filter over $P(U)$. Moreover, if $C=C_G$ then $G_C=G$.
\end{lemma}
\textit{Proof}: 
Let $p\in G_C$ and $p'\leq p$, then
$$t(p)=C\cap (\max(t(p))+1), \ C\setminus (\max(t(p))+1)\subseteq
A(p).$$ 
To see that $p'\in G_C$,
$$C\cap (\max (t(p'))+1)=C\cap(\max(t(p))+1)\cap(\max(t(p'))+1)=t(p)\cap(\max(t(p'))+1)=t(p').$$
Moreover,
$$C\setminus (\max(t(p'))+1)=[C\setminus( \max(t(p))+1)]\cup [t(p)\setminus(\max(t(p'))+1)]\subseteq A(p)\cup A(p')=A(p'),$$
hence $p'\in G_C$.
Suppose $p,q\in G_C$, without loss of generality assume that $n(p)\geq n(q)$.
Since both $t(p),t(q)$ are initial segments of $C$, $t(p)$ is an end-extension of $t(q)$.
Since $$C\setminus( \max(t(p))+1)\subseteq C\setminus( \max(t(q))+1)\subseteq A(q),$$ it follows that $$C\setminus (\max(t(p))+1)\subseteq A(p)\cap A(q)\text{ and } p^{\frown}\langle A(p)\cap A(q)\rangle\in G_C$$ 
Moreover, $p,q\leq p^{\frown}\langle A(p)\cap A(q)\rangle.$ This conclude the proof that $G_C$ is a filter.
Now suppose that $C=C_G$ is Prikry-generic. Since generic filters are maximal among filter,
it suffices to prove that $G\subseteq G_C$. Let $p\in G$, in the proof of lemma \ref{Genericseq} $(3)$ we have seen that $$C_G\setminus(\max(t(p))+1)\subseteq A(p).$$
Obviously, $t(p)\subseteq C_G\cap(\max(t(p))+1)$. For the other inclusion, let $\alpha\in C_G\cap(\max(t(p))+1)$, there is $r\in G$ such that $\alpha\in r$ and $r\geq p$. Then $\alpha\in t(r)\cap\max(t(p))+1=t(p)$, thus $ t(p)=C_G\cap\max(t(p))+1$ and $p\in G_C$.
$\square$

Now we turn to the preservation of cardinals. In \cite{Devlin}, Devlin proved that the filters $U$ for which $P(U)$ does not collapse cardinals are exactly the shrinking ultrafilters.  We include here the proof of one direction and refer the reader to \cite{Devlin} for the other.

\begin{proposition}
If $U$ is shrinking, then $P(U)$ preserves cardinals.
\end{proposition}
\textit{Proof:} Since $P(U)$ is $\kappa^+$-c.c. and $
\kappa$ is a limit cardinal, it suffices to show that successor cardinals $<\kappa$ are preserved. Let $\lambda^+<\kappa$ be a successor cardinal. Let $f:\lambda\rightarrow \lambda^+\in V[G]$ be any increasing function and $\dot{f}$ be a $P(U)$-name such that $$\Vdash_{P(U)} \dot{f}:\lambda\rightarrow\lambda^+ \text{  is  increasing}.$$ Let $p\in P(U)$, $\xi<\lambda$
and $\vec{\alpha}\in[A(p)]^{<\omega}$. There is at most one $\nu$ such that for some $p^{\frown}\vec{\alpha}\leq^* q$,
$q\Vdash \dot{f}(\xi)=\nu$. To see this, note that any $p^{\frown}\vec{\alpha}\leq^* q,q'$ are compatible and cannot force contradictory information. Define
$$E_{\vec{\alpha}}=\{\nu<\lambda^+\mid \exists \xi<\lambda\exists q\geq^* p^{\frown}\vec{\alpha}, \ q\Vdash \dot{f}(\xi)=\nu\}.$$
Thus for any $\vec{\alpha}\in [A(p)]^{<\omega}$, $|E_{\vec{\alpha}}|<\lambda^+$. In $V$, $\lambda^+$ is regular the $\gamma_{\vec{\alpha}}:=\sup(E_{\vec{\alpha}})<\lambda^+$. Define $$F:[A(p)]^{<\omega}\rightarrow \lambda^+
, \ F(\vec{\alpha})=\gamma_{\vec{\alpha}}.$$ By  the shrinking property of $U$, there is $A^*\in U$ such that  $|F''[A^*]^{<\omega
}|<\lambda^+$. Since $\lambda^+$ is regular in $V$, $$\gamma^*:=\sup(F(\vec{\alpha})\mid \vec{\alpha}\in[A^*]^{<\omega})<\lambda^+.$$ 
Let us prove that $p^*=p^{\frown}\langle A^*\rangle\in P(U)$ forces that  $Im(\dot{f})\subseteq \gamma^*$.
Indeed, if  $ p^{*\frown}\vec{\alpha}\leq^*q$  and $q\Vdash \dot{f}(\xi)=\nu$,
then $\nu\in E_{\vec{\alpha}}$. By definition, $\nu<\gamma_{\vec{\alpha}}$ and since $\vec{\alpha}\in [A^*]^{<\omega}$, $\nu_{\vec{\alpha}}<\gamma^*$.
By density, we can find such a condition $p^*\in G$ and conclude that $f=(\dot{f})_G$ is bounded.
$\square$
\begin{corollary} For any  normal or Rowbottom ultrafilter $U$,
$P(U)$ preserve cardinals.
\end{corollary}

Prikry's original forcing uses a normal ultrafilter and has the other crucial property that $P(U)$ does not add new bounded subsets to $\kappa$. Already in \cite{Prikry}, Prikry noticed that a Rowbottom ultrafilter  suffices to conclude the preservation of cardinals without adding bounded subsets to $\kappa$. What is the exact requirement of $U$ so that $P(U)$ will not add new bounded subsets to $\kappa$? 
First note that $U$ must be an ultrafilter.
\begin{proposition}\label{newreal}
If $U$ contains all final segments and is not an ultrafilter, then $U$ adds a new real. 
\end{proposition}
\textit{Proof}: Let $Z\subseteq\kappa$ such that $Z,\kappa\setminus Z\notin U$. Let $G$ be generic and consider the set $$X=\{n<\omega\mid C_G(n)\in Z\}\in V[G].$$ We claim that $X\notin V$. Otherwise, Let $\dot{X}$ be a name for $X$ and find $p\in G$ such that $p\Vdash \dot{X}=X$. It follows that $$\forall  i\in n(p)\cap X, \ t_i(p)\in Z.$$ Since $A(p)\in U$,  $$Z\cap A(p),(\kappa\setminus Z)\cap A(p)\neq\emptyset.$$ Let $x\in Z\cap A(p)$  and $y\in (\kappa\setminus Z)$. If $n(p)+1\in X$, consider the condition $p^{\frown}\langle y\rangle\geq p$. Then
$$p^{\frown}\langle y\rangle\Vdash n(p)+1\in\dot{X}\wedge n(p)+1\notin \dot{X},$$
contradiction. The case $n(p)+1\in \kappa\setminus X$ is similar. 
$\square$

 Now we turn to the Prikry condition, which will be proven for Rowbottom ultrafilters.
Devlin and Paris established that the Rowbottom property is also necessary \cite{Devlin}.
\begin{theorem}
Suppose $U$ is a Rowbottom ultrafilter, let $\sigma$ be a formula in the forcing language and $p\in P(U)$. Then there exists $p\leq^* p^*$ such that $p^* \ || \sigma$ i.e. $p\Vdash\sigma\vee p\Vdash\neg\sigma$.
\end{theorem}
 \textit{Proof:} Define a function $F:[A(p)]^{<\omega}\rightarrow \{0,1,2\}$ by 
$$ 
F(\vec{\alpha})=\begin{cases}
        0, & \exists A', \ p^{\frown}\langle\vec{\alpha},A'\rangle\Vdash \sigma\\
        1, & \exists A', \ p^{\frown}\langle\vec{\alpha},A'\rangle\Vdash \neg\sigma\\
        2, & \text{Otherwise } \end{cases}.$$
The function $F$ is well defined since for every $t\in[\kappa]^{<\omega}$ and $A_1,A_2\in U$, the conditions $\langle t,A_1\rangle,\langle t,A_2\rangle$ are compatible and cannot force contradictory information such as $\sigma$ and $\neg\sigma$. Since $U$ is  Rowbottom, there is $A^*\in U$ such that for all $n<\omega$, $F\restriction [A^*]^n$ is constant.
Let $p^*=p^{\frown}\langle A^*\rangle$, to see that $p^* || \sigma$, suppose otherwise, then there are $s_1,s_2\in[A^*]^{<\omega}$ and $X_1,X_2\in P(A^*)\cap U$ such that $$p^{*\frown}\langle s_1,X_1\rangle\Vdash \sigma, \ \ \ p^{*\frown}\langle s_2,X_2\rangle\Vdash\neg\sigma.$$ Extend $s_1,s_2$ to $s^*_1,s^*_2$ if necessary, so that $|s^*_1|=|s^*_2|=n^*$. The extended conditions $p^{*\frown}\langle s^*_i,X_i\rangle$ decide $\sigma$ the same way $p^{*\frown}\langle s_i,X_i\rangle$ did. We have reached a contradiction since $s^*_1,s^*_2\in [A']^{n^*}$ and $F(s^*_1)=0, \ F(s^*_2)=1$. 
$\square$
\begin{corollary}
For any normal or Rowbottom ultrafilter $U$, if $\nu<\kappa$ then $P^{V[G]}(\nu)=P^{V}(\nu)$.
\end{corollary}

The next lemma is known as the strong Prikry property.
\begin{lemma}\label{strongPrikry}
Let $U$ be a normal ultrafilter and $D\subseteq P(U)$ dense and open. Then for every $p\in P(U)$, there is $p\leq^*p^*$ and a natural number $n<\omega$ such that $$\forall\vec{\alpha}\in [A(p^*)]^n, \ p^{*\frown}\vec{\alpha}\in D.$$
\end{lemma}
\textit{Proof}: Define $F:[A(p)]^{<\omega}\rightarrow \{0,1\}$ by $$F(\vec{\alpha})=0\Longleftrightarrow  \ \exists A\in U, \ p^{\frown}\langle\vec{\alpha},A\rangle\in D.$$
By proposition \ref{BasicProp}, $U$ is Rowbottom. Therefore, there is $A'\in P(A(p))\cap U$ such that for every $n<\omega$, $F\restriction[A']^n$ is constant. If $F(\vec{\alpha})=0$, let $A_{\vec{\alpha}}\in U$ be a witnessing set for it. Define 
$$A^*=A'\cap(\triangle_{\vec{\alpha}\in[A']^{<\omega}}A_{\vec{\alpha}}):=\{\nu\in A'\mid \forall\vec{\alpha}\in[A']^{<\omega}(\max(\vec{\alpha})<\nu\rightarrow\nu\in A_{\vec{\alpha}})\}.$$
By normality, $A^*\in U$. Since $D$ is dense, there is $n$ and $\vec{\alpha}\in [A^*]^n$ such that for some $A\in U$
$p^{\frown}\langle\vec{\alpha},A\rangle\in D$. It follows that $F(\vec{\alpha})=0$.
Let us claim that $p^*=p^{\frown}\langle A^*\rangle$ and $n$ are as wanted. 
Fix $\vec{\beta}\in[A^{*}]^{n}$, then
$F(\vec{\beta})=F(\vec{\alpha})=0$. By the definition of $F$,
$$F(\vec{\beta})=0\Rightarrow \exists B \  p^{\frown}\langle\vec{\beta},B\rangle\in D\Rightarrow p^{\frown}\langle\vec{\beta},A_{\vec{\beta}}\rangle\in D.$$
Note that $A^*\setminus(\max(\vec{\beta})+1)\subseteq A_{\vec{\beta}}$, hence $p^{\frown}\langle\vec{\beta},A_{\vec{\beta}}\rangle\leq p^{*\frown}\langle \vec{\beta}\rangle$. Finally, since $D$ is open it follows that $p^{*\frown}\langle\vec{\beta}\rangle\in D$.
$\square$

The following theorem is due to Mathias \cite{Mathias}. 
\begin{theorem}\label{Mathias}
 Let $U$ be normal. Suppose  $C\in[\kappa]^{\omega}$ is a sequence such that:
 \begin{enumerate}
 \item $\sup(C)=\kappa$.
 \item For every $A\in U$, there is $n<\omega$ such that for every $m\geq n$, $C(m)\in A$.
 \end{enumerate}
Then $C$ is Prikry-generic for $U$.
\end{theorem}
\textit{Proof}: Let us show that $G=G_C$ is generic. By lemma \ref{Filter} it is a filter. Let $D$ be a dense open subset of $P(U)$. For every $s\in[\kappa]^{<\omega}$, apply the strong Prikry property to the condition $\langle s,\kappa\setminus \max(s)+1\rangle$ and $D$,  find $A_s\in U$ and $n_s<\omega$ such that $$\forall\vec{\alpha}\in [A_s]^{n_s}, \ \langle s^{\frown}\vec{\alpha},A_s\setminus(\max(\vec{\alpha})+1)\rangle\in D.$$ Let $$A^*=\triangle_{s\in[\kappa]^{<\omega}} A_s\in U.$$ By requirement (2) of the sequence $C$, there is $m<\omega$ such that for every $n\geq m$, $C(n)\in A^*$. 
Consider the condition
$$p=\langle C\restriction m+1,A^*\rangle$$ and let $n^*=n_{C\restriction m+1}$. By the definition of $G$,
$p^{\frown}\langle C(m+1),...,C(n^*)\rangle\in G$. Also $$A^*\setminus (C(m)+1)\subseteq A_{C\restriction m+1},$$ therefore $$\langle C(m+1),...,C(n^*)\rangle\in [A_{C\restriction m+1}]^{n^*}.$$ Since $D$ is open,  $p^{\frown}\langle C(m+1),...,C(n^*)\rangle\in D\cap G$.
$\square$
\begin{corollary}\label{othergeneric}
Let $C,D$ be Prikry-generic sequences for some Rowbottom ultrafilter $U$. Then infinite subsequences of $C$ and the sequence $C\cup D$ are again Prikry-generic for $U$.
 \end{corollary}
\textit{Proof}: These sequences satisfy the Mathias criteria.
$\square$

The following corollary is due to Solovay,  Dehornoy in \cite{Dehornoy} and Bukovsky in \cite{Bukovsky} independently. We follow the notations of chapter 19 in \cite{Jech}.
\begin{corollary}
Let $U$ be a normal measure on $\kappa$ and denote by $$\langle Ult_{U}^{(n)},\kappa^{(n)},U^{(n)},i_{n,m}\mid n<m\leq\omega\rangle,$$ the system of the first $\omega$ iterations of the iterated ultrapower by $U$, then:
\begin{enumerate}
\item $\langle \kappa^{(n)}\mid n<\omega\rangle$ is Prikry-generic for $U^{(\omega)}$ over the model $Ult^{(\omega)}_U$.
\item $Ult^{(\omega)}_U[\langle \kappa^{(n)}\mid n<\omega\rangle]=\bigcap_{n<\omega}Ult^{(n)}_U$.
\end{enumerate}
\end{corollary}
\textit{Proof}: For (2), see \cite{Dehornoy} or \cite{Bukovsky}.
For (1), we prove that the Mathias criteria holds for the critical sequence $\langle\kappa^{(n)}\mid n<\omega\rangle$. 
Let $X\in U^{(\omega)}$, by the definition of $Ult^{(\omega)}_U$ as the direct limit of the system $$\langle Ult^{(n)}_U,i_{n,m}\mid n<m<\omega\rangle,$$ there is $n<\omega$ and $Y\in Ult^{(n)}_U$ such that $i_{n,\omega}(Y)=X$.
Fix any $n\leq m<\omega$ and denote $Z=i_{n,m}(Y)$. By normality of $U^{(m)}$, $\kappa^{(m)}\in i_{m,m+1}(Z)$,
by elementarity, $$i_{m+1,\omega}(\kappa^{(m)})\in i_{m+1,\omega}(i_{m,m+1}(Z))=i_{n,\omega}(Y)=X.$$ Moreover, $\kappa^{(m)}<\kappa^{(m+1)}=crit(i_{m+1,\omega})$, then $\kappa^{(m)}\in X$.
$\square$
\begin{corollary}

If $U$ is Rowbottom, then $P(U)$ has the Mathias characterization.
\end{corollary}
\textit{Proof}: By proposition \ref{Rownor}, there is a normal ultrafilter $W$ such that $U\equiv_{R-K}W$. Let $f:\kappa\rightarrow\kappa$ be a bijection such that $f_*(U)=W$ and $f^{-1}_*(W)=U$. Assume that $C$ is such that for every $A\in U$, $C\setminus A$ is finite. Then for every $A\in U$, $f[C\setminus A]$ is finite. Since $f_*(U)=W$, for every $X\in W$,  $f[C]\setminus X$ is finite and therefore by theorem \ref{Mathias}, $f[C]$ is $W$-generic. By proposition \ref{Genericseq}, $C=f^{-1}[f[C]]$ is also generic.
 $\square$
 \begin{theorem}
 If $P(U)$ satisfies the Mathias characterization, then $U$ is Rowbottom.
 \end{theorem}
 \textit{Proof}: Suppose $U$ is not Rowbottom. By \ref{Rownor}, there is $f:[\kappa]^{2}\rightarrow 2$ with no homogeneous set from $U$. Let $G$ be generic for $P(U)$. In $V[G]$, define $$G:[\omega]^2\rightarrow 2, \ G(n,m)=f(C_G(n),C_G(m)).$$ By Ramsey's theorem (see \cite{Jech}, theorem 9.1), there is $H\subseteq\omega$ such that $|H|=\omega$, homogeneous for $G$. Consider $C'=\{C_G(n)\mid n\in H\}$, then $C'$ is homogeneous for $f$. Moreover, For every $A\in U$, $$C'\setminus A\subseteq C_G\setminus A.$$ 
 By \ref{Genericseq}, $C_G\setminus A$ is finite, hence $C'$ satisfies the Mathias criteria. Let us show that $C'$ is not generic. To see this, we will prove that generic sequences cannot be homogeneous for $f$. Toward a contradiction, assume that
 $$p\Vdash\dot{C} \text{  is  homogeneous  for } f,$$
 where $\dot{C}$ is a name for the Prikry-generic sequence. Since $A(p)$ is not homogeneous for $f$,  we can choose $\alpha,\beta,\gamma,\delta\in A(p)$ such that $f(\alpha,\beta)\neq f(\gamma,\delta)$. Assume $\alpha<\beta<\gamma<\delta$, other possibilities of ordering these four ordinals are treated similarly. Extend $p$ to the condition $p^*=p^{\frown}\langle\alpha,\beta,\gamma,\delta\rangle$, then $$p^*\Vdash \alpha,\beta,\gamma,\delta\in\dot{C}.$$ This is a contradiction since $p^*$ also forces that $\dot{C}$ is homogeneous.   It follows that $C'$ is a counterexample for the Mathias characterization.
 $\square$
 
 The next lemma is a topological separation property which to the best of our knowledge is unknown. 
 \begin{lemma}
 If $U$ is a normal ultrafilter over $\kappa$, $X\in V[G]$ is a set such that $X\cap C_G=\emptyset$ and $|X|<\kappa$, then there is $A\in U$ such that $A\cap X=\emptyset$.
 \end{lemma}
 \textit{Proof}:
 It suffices to show that if $\dot{x}$ is a name such that $p\Vdash \dot{x}\notin \dot{C_G}$, then there if $p\leq^*p^*$ and $A\in U$ such that $p^*\Vdash \dot{x}\notin A$. Once  we have established that, then the argument for a general set $X$ is simple, let  $\langle\dot{x_i}\mid i<|X|\rangle$ be a name for an enumeration of $X$. 
 Construct a $\leq^*$-increasing sequence $\langle p_i\mid i<|X|\rangle$ and sets $A_i\in U$ such that $$p_i\Vdash \dot{x_i}\notin A_i.$$ Since $U$ is $\kappa$-complete, find $p^*$ and $A^*$   such that $p_i\leq^*p^*$ and $p^*\Vdash\dot{X}\cap A^*=\emptyset$. 
 By density, such $p^*$ can be found in $G$ and so $A^*\cap X=\emptyset$. Assume $p\Vdash \dot{x}\notin \dot{C_g}$, define $$f:[A(p)]^{<\omega}\rightarrow \{0,1\}, \  \ f(\vec{\alpha})=0\Leftrightarrow\exists C \ p^{\frown}\langle \vec{\alpha},C\rangle || \dot{x}.$$ Let $A_1\subseteq A(p)$ be homogeneous and let $n<\omega$ be minimal such that $f\restriction[A_1]^n\equiv0$. For $\vec{\alpha}\in[A_1]^n$, fix $E_{\vec{\alpha}},x_{\vec{\alpha}}$ such that $p^{\frown}\langle \vec{\alpha},E_{\vec{\alpha}}\rangle\Vdash \dot{x}=x_{\vec{\alpha}}$. Define
 $$A_2=\Delta_{\vec{\alpha}\in[A_1]^n}E_{\vec{\alpha}}.$$
 It follows that for every $\vec{\alpha}\in [A_2]^{<\omega}$,  $x_{\vec{\alpha}}\notin A_2\setminus (\max(\vec{\alpha})+1)$. Otherwise extend $p$ to $p^{\frown}\langle \vec{\alpha}, x_{\vec{\alpha}}, A_2\rangle$, 
 this condition forces $x_{\vec{\alpha}}\in \dot{C_G}\ \wedge \  x_{\vec{\alpha}}=\dot{x}$ but also $\dot{x}\notin C_G$, contradiction. Shrink $A_2$ to $A_3\in U$, so that for every $\vec{\alpha}\in[A_3]^n$, either $x_{\vec{\alpha}}>\max(\vec{\alpha})$ or $x_{\vec{\alpha}}<\max(\vec{\alpha})$, equality cannot hold since $$p^{\frown}\vec{\alpha}\Vdash\max(\vec{\alpha})\in \dot{C_G}.$$
 Assume that $\max(\vec{\alpha})>x_{\vec{\alpha}}$, shrink $A_3$ to $A_4\in U$ and find $i<n$ such that for every $\vec{\alpha}\in[A_4]^n$, $\alpha_i<x_{\vec{\alpha}}<\alpha_{i+1}$ (given that $\alpha_0=0$). Fix $\langle\alpha_1,...,\alpha_i\rangle$, then the function $$\langle\alpha_{i+1},...,\alpha_n\rangle\mapsto x_{\vec{\alpha}}$$ is shrinking. By normality of $U$, there is $A_{\alpha_1,...,\alpha_i}\in U$ and $y_{\alpha_1,...,\alpha_i}$ so that $$\forall \vec{\gamma}\in[A_{\alpha_1,...,\alpha_i}]^{n-i}, \ x_{\alpha_1,...,\alpha_i,\vec{\gamma}}=y_{\alpha_1,...,\alpha_i}.$$
 Let $$A_5=\triangle_{\vec{\beta}\in [A_4]^i}A_{\vec{\beta}}\in U.$$ It follows that for every $\vec{\beta}\in[A_5]^i$, $p^{\frown}\vec{\beta}\Vdash \dot{x}=y_{\vec{\beta}}$. By the same argument as before, $y_{\vec{\beta}}\notin A_5\setminus (\max(\vec{\beta})+1)$, but also $y_{\vec{\beta}}>\max(\vec{\beta})$, hence $y_{\vec{\beta}}\notin A_5$.
 We claim that $p^{\frown}\langle A_5\rangle\Vdash \dot{x}\notin A_5$. Otherwise, there is $q\geq p^{\frown}\langle A_5\rangle$ forcing $\dot{x}\in A_5$. Extend $q$ if necessary so that for some
 $\vec{\beta}\in[A_5]^{i}$,  $q\geq p^{\frown}\langle\vec{\beta},A_5\rangle$. It follows that $$p^{\frown}\langle \vec{\beta},A_5\rangle\Vdash \dot{x}=y_{\vec{\beta}}.$$
 But then $y_{\vec{\beta}}\in A_5$, contradiction.
 $\square$
 \begin{corollary}
 If $c\in V[d]$, where $c,d$ are Prikry-generic for the same $U$, then $c\setminus d$ is finite. 
 \end{corollary}
 \textit{Proof}: Otherwise use the previous lemma and the Mathis characterization to obtain a contradiction.
 $\square$
 \begin{corollary}
 For any two Prikry-generic sequences, $c$ and $d$ for the same $U$,  $$V[c]=V[d] \ \Leftrightarrow \ |c\Delta d|<\omega.$$
 \end{corollary}

 The Mathias criteria and the previous corollaries suggests that Prikry generic extensions with a Rowbottom ultrafilter are in some sense rigid. This leads to the natural question:
 \begin{center}
 \textit{What are the intermediate ZFC models $V\subseteq M\subseteq M[G]$?}
 \end{center}
 In \cite{PrikryCase}, Gitik, Koepke and Kanovei answered this question:
 \begin{theorem}\label{ParikryCaseInter}
 Let $U$ be a normal ultrafilter and let $G$ be $P(U)$-generic. Then for every ZFC model $V\subseteq M\subseteq V[G]$, there is $C'\subseteq C_G$ such that $V[C']=M$. 
 \end{theorem}
 
 By corollary \ref{othergeneric}, it follows that:
 \begin{corollary}
 Every intermediate model of a Prikry generic extension is again Prikry generic.
 \end{corollary}
 
\section{Tree-Prikry Forcing}

The Tree-Prikry forcing is an alternative forcing to Prikry forcing that uses ultrafilters and filters that need not be Rowbottom. 
\begin{definition}\label{TreeOfFilters}
A \textit{tree of filters} is a sequence $\vec{U}=\langle U_a\mid a\in[\kappa]^{<\omega}\rangle$, such that for every $a\in[\kappa]^{<\omega}$, $U_a$ is a filter over $\kappa$.
\end{definition} 
\begin{definition}\label{FatTree}
Let $\kappa$ be a measurable cardinal and let $\vec{U}=\langle U_a\mid a\in[\kappa]^{<\omega}\rangle$ be a tree of filters. A tree $\langle T,\leq_T\rangle$ is a \textit{$\vec{U}$-fat tree with trunk $t$} if:
\begin{enumerate}
\item $T\subseteq[\kappa]^{<\omega}$.
\item $\forall\vec{\alpha},\vec{\beta}\in T$, \  $\vec{\alpha}\leq_T\vec{\beta}\Leftrightarrow\vec{\beta}\cap\max(\vec{\alpha})+1=\vec{\alpha}$.
\item $T$ is $\leq_T$-downward closed.
\item For every $t\in T$, either $t<_T t^*$ or  $\succ_T(t):=\{\alpha<\kappa\mid t^{\frown}\alpha\in T\}\in U_t$.
\end{enumerate}
\end{definition}

We use other standard notation of trees:
\begin{enumerate}
\item For $t\in T$ such that $t^*\leq t$, $ht_T(t)=|t|-|t^*|$.
\item $\Lev_n(T)=\{t\in T\mid ht_T(t)=n\}$.
\item For $A\subseteq\omega$,  $T\restriction A=\bigcup_{a\in A}\Lev_a(T)$.
\item For $t\in T$, $(T)_t=\{s\in T\mid s\geq_T t\vee s<_T t\}$. If $t^*$ is the trunk of $T$ and $t^*\leq t$, then $(T)_t$ is a $\vec{U}$-fat tree with trunk $t$.
\item $mb(T)=\{b\in[\kappa]^{\omega}\mid \forall n<\omega \  b\restriction n \in T\}$ is the set of maximal branches of $T$.
\end{enumerate}
\begin{proposition}\label{PropTree}
Let $T$ be a $\vec{U}$-fat tree with trunk $t$:
\begin{enumerate}
\item Let $\langle A_a\mid a\in[\kappa]^{<\omega}\rangle$ such that $A_a\in U_a$. Then for every $\vec{U}$-fat tree $T$, there is a $\vec{U}$-fat tree $T^*\subseteq T$, such that for every $a\in T^*$, $\succ_{T^*}(a)\subseteq A_a$.
\item Assume that all the filters in $\vec{U}$ are $\lambda$-complete. Let $\beta<\lambda$ and $\langle T_\alpha\mid \alpha<\beta\rangle$ be a sequence of $\vec{U}$-fat trees with the same trunk $t$. Then $\cap_{\alpha<\lambda}T_\alpha$ is also a $\vec{U}$-fat tree with trunk $t$.

\item Let $T$ be a $\vec{U}$-fat tree with trunk $t$. Suppose  $\langle T(s)\mid s\in T \wedge s>_T t\rangle$ is a tree of $\vec{U}$-fat trees such that $T(s)\subseteq (T)_s$ has trunk $s$. Then we can amalgamate them to a single $\vec{U}$-fat tree i.e. there is a $\vec{U}$-fat tree $T^*$ with trunk $t$, such that for every $s\in T^*$, $(T^*)_s\subseteq T(s)$.

\end{enumerate}
\end{proposition}
\textit{Proof}: For $(1)$, we define $T^*$ by induction on $\Lev_n(T^*)$. $T^*$ has the same trunk as $T$. 
Assume $s\in \Lev_n(T^*)$ and define $$\succ_{T^*}(s)=A_s\cap \succ_{T}(s)\in U_s.$$ 
For $(2)$,  set $T'=\cap_{\alpha<\lambda}T_{\alpha}$ and let us verify that definition \ref{FatTree} holds for $T'$.
For $s\in T'$ either $t$ is an end-extension of $s$ and $s\leq_{T'}t$, or $t\leq_{T_\alpha}s$ for every $\alpha<\beta$, implying that $\succ_{T_{\alpha}}(s)\in U_s$. 
By $\lambda$-completeness of $U_s$, $\succ_{T'}(s)=\cap_{\alpha<\lambda}\succ_{T_\alpha}(s)\in U_s$. Finally for $(3)$,
let us define inductively trees with the same trunk $t$:
$$...\subseteq T^{(n)}\subseteq T^{(n-1)}\subseteq...\subseteq T^{(0)}$$
such that fro every $n<\omega$, $T^{(n+1)}\restriction (n+1)=T^{(n)}\restriction (n+1)$. We start with  $T(t)=T^{(0)}$. Assume we have defined $T^{(n)}$ and let 
$$T^{(n+1)}=(T^{(n)}\restriction n+1)\cup\bigcup_{s\in \Lev_n(T^{(n)}) }(T(s)\cap (T^{(n)})_s).$$
Namely, $T^{(n+1)}$ is the same as $T^{(n)}$ up to the $n$th level and for every $s\in \Lev_n(T^{(n)})$, we shrink the tree $(T^{(n)})_s$ to $(T^{(n)})_s\cap T(s)$. By $(2)$, $T(s)\cap (T^{(n)})_s$ is a $\vec{U}$-fat tree with trunk $s$, it follows that definition \ref{FatTree} is satisfied for the tree $T^{(n+1)}$.
Let $T^*=\cap_{n<\omega}T^{(n)}$, note that the filters involved need not be even $\sigma$-complete,
but $T^{(n+1)}\restriction (n+1)=T^{(n)}\restriction (n+1)$ and the intersection at each level is in fact of finitely many sets.
It follows that $T^*$ is a $\vec{U}$-fat tree with trunk $t$ and for every $s\in T^*$, $$(T^*)_s\subseteq (T^{(ht(s))})_s\subseteq T(s).$$
$\square$
\begin{definition}

Let $\vec{U}=\langle U_a\mid a\in[\kappa]^{<\omega}\rangle$ be a tree of filters over $\kappa$. 
The underlining set of $P_T(\vec{U})$ is the set of all $\langle t_1,...,t_n,T\rangle$, where $ \langle t_1,...,t_n\rangle$ is the trunk of the $\vec{U}$-fat tree $T$. For $p=\langle t_1,..., t_n ,T\rangle, \  q=\langle s_1,...,s_m,S\rangle$, define $p\leq q$ iff:
\begin{enumerate}
\item $n\leq m$.
\item $\langle s_1,..., s_m\rangle\in T$ (in particular $\langle s_1,..., s_m\rangle$ is an end-extension of $\langle t_1,...,t_n\rangle$).
\item $S\subseteq T$.
\end{enumerate}

Let $p\leq^* q$ iff $p\leq q\wedge n=m$.
\begin{enumerate}
\item  For $\vec{W}=\langle W_\alpha \mid \alpha<\kappa\rangle$, define $P_T(\vec{W})=P_T(\vec{U})$, where $U_a=W_{\max(a)}$ for all $a\in[\kappa]^{<\omega}$.
If the $W_\alpha$'s are distinct normal measures on $\kappa$ this is in fact the minimal Prikry forcing appearing in \cite{MinimalPrikry}.
\item For $\vec{W}=\langle W_n\mid n<\omega\rangle$
define $P_T(\vec{W})=P_T(\vec{U})$, where $U_a=U_{|a|}$ for all $ a\in [\kappa]^{<\omega}$.
\item For a single filter $W$, define $P_T(W)=P_T(\vec{U})$ where $U_a=W$ for all $a\in[\kappa]^{<\omega}$.
\end{enumerate}
\end{definition}
We will refer to those forcings by \textit{Tree-Prikry with a tree of filters}, \textit{Tree-Prikry with a $\kappa$-sequence of filters}, \textit{Tree-Prikry with an $\omega$-sequence of filters} and \textit{Tree-Prikry with a single filter}.

The following notations are useful for a condition $p=\langle t_1,...,t_n,A\rangle\in P_T(\vec{U})$:
\begin{enumerate}
\item  $T=T(p)$.
\item $t(p)=\langle t_1,...,t_n\rangle$.
\item $n(p)=n$.
\item For $i\leq n$, $t_i(p)=t_i$.
\end{enumerate}
If $t(p)^{\frown}\vec{\alpha}\in\Lev_m(T(p))$ and $T'\subseteq T(p)$ is a $\vec{U}$-fat tree with trunk $t(p)^{\frown}\vec{\alpha}$, define:
\begin{enumerate}
\item $p^{\frown}\langle\vec{\alpha},T'\rangle=p^{\frown}\langle t(p)^{\frown}\vec{\alpha},T'\rangle:=\langle t(p)^{\frown}\vec{\alpha},T'\rangle$.
\item $p^{\frown}\vec{\alpha}=p^{\frown}(t(p)^{\frown}\vec{\alpha}):= p^{\frown}\langle\vec{\alpha},(T(p))_{\vec{\alpha}}\rangle$. 

\end{enumerate}  \begin{proposition}
If $U$ is a normal ultrafilter over $\kappa$, then $P_T(U)$ and $P(U)$ are forcing equivalent.
\end{proposition}
\textit{Proof}: Define $\pi:P(U)\rightarrow P_T(U)$ by 
$$\pi(p)=\langle t(p), T(\pi(p))\rangle.$$
Where $T(\pi(p))$ is a tree with trunk $t(p)$ and for every $t\in T(p)$ such that $t(p)\leq t$, $\succ_{T(p)}(t)=A\setminus(\max(t)+1)$. It follows that $\pi(p)\in P_T(U)$. To see that $Im(\pi)$ is dense in $P_T(U)$, fix $p\in P_T(U)$,  define the set $$A=\triangle_{t(p)\leq t\in T(p)}\succ_{T(p)}(t).$$
By normality of $U$, $A\in U$. Consider the condition $q=\langle t(p),A\rangle$, then $q\in P(U)$. Moreover $T(\pi(q))\subseteq T(p)$. Indeed, let $t(p)^{\frown}\langle\alpha_1,...,\alpha_n\rangle\in T(\pi(q))$, then for every $1\leq i< n$, $$\alpha_{i+1}\in A\setminus(\alpha_{i}+1)\subseteq \succ_{T(p)}(t(p)^{\frown}\langle\alpha_1,...,\alpha_i\rangle).$$ In particular, $t(p)^{\frown}\langle\alpha_1,...,\alpha_n\rangle\in T(p)$. It follows that $\pi(q)\geq p$ and so $Im(\pi)$ is dense in $P_T(U)$.
 $\pi$ is order preserving, let $p,q\in P(U)$ such that $p\leq q$, then:
 \begin{enumerate}
     \item $t(q)$ is an end-extension of $t(p)$.
     \item $t(q)\setminus t(p)\in[A(p)]^{<\omega}$.
     \item $A(q)\subseteq A(p)$.
 \end{enumerate} 
 It follows that  $T(\pi(q))\subseteq T(\pi(p))$ and $t(q)\in T(\pi(p))$, therefore $\pi(p)\leq \pi(q)$.
 Assume that $p,q\in P(U)$ and $\pi(p),\pi(q)$ have a common extension $r\in P_T(U)$. By definition of the order, $t(r)$ is an end-extension of $t(p),t(p)$. Moreover, since $t(r)\in T(\pi(p))\cap T(\pi(q))$, $t(r)\setminus t(q)\in [A(q)]^{<\omega}$ and $t(r)\setminus t(p)\in [A(p)]^{<\omega}$ and the condition $$\langle t(r),A(p)\cap A(q)\setminus(\max(t(r))+1)\rangle\in P(U)$$
 is a common extension of $p,q$.
$\square$
\begin{lemma} For any tree of filters $\vec{U}=\langle U_a\mid a\in [\kappa]^{<\omega}\rangle$:
\begin{enumerate}
\item $P_T(\vec{U})$ is $\kappa^+$-Knaster. 
\item If all the filters in $\vec{U}$ are $\lambda$-complete, then the order $\leq^*$ is $\lambda$-closed.
\end{enumerate}
\end{lemma}
\textit{Proof}:
The proof is completely analogous to \ref{PrikryProp}.
$\square$

To avoid similar pathologies to those of example \ref{pathologies}, assume that all the filters in $\vec{U}$ contains all the final segments. The next theorem is a variation the Rowbottom property for a tree of filters.
\begin{theorem} \label{RowTree}
Let $T$ be a $\vec{U}$-fat tree and $F:T\rightarrow X$ a function such that for every $a\in[\kappa]^{<\omega}$, $U_a$ is a $|X|^+$-complete ultrafilter. Then there exists a $\vec{U}$-fat tree $T'\subseteq T$ with the same trunk, such that for any $n<\omega$, $F\restriction \Lev_n(T')$ is constant.
\end{theorem}
\textit{Proof}: 
By induction on $n<\omega$, we will prove that for every such $T$ and $F$, we can find $T'\subseteq T$ such that for every $m\leq n$, $F\restriction \Lev_m(T')$ is constant.  
For $n=1$, let $t^*$ be the trunk of $T$ and consider $F\restriction \succ_T(t^*)$.
Since $U_{t^*}$ is a $|X|^+$-complete ultrafilter, there is $A\in U_{t^*}$ such that $F\restriction A$ is constant.
Shrink the set $\succ_T(t^*)$ to $A$ and denote the new tree by $T^{(1)}$, this tree is as wanted. 
Assume the theorem hold up to $n$ and let $T$ and $F$ be as before. 
By the induction hypothesis, there is $T^{(n)}\subseteq T$ such that for every $m\leq n$, $F\restriction \Lev_m(T^{(n)})$ is constant.
For $t\in \Lev_n(T^{(n)})$, define $$F_t:\succ_{T^{(n)}}(t)\rightarrow X, \ F_t(\alpha)=F(t^{\frown}\langle\alpha\rangle).$$
There is $A_t\in U_t$ such that $F_t\restriction A_t$ is constant with value $\gamma_t$. Define the tree $T'$,
$$T'=(T^{(n)}\restriction(n+1))\cup\bigcup_{t\in \Lev_n(T^{(n)}),\ \alpha\in A_t} (T^{(n)})_{t^{\frown}\alpha}.$$
Namely, $T'$ is the same as $T^{(n)}$ up to the $n$th level, then we shrink the tree so that the sets $A_t$ constitute $\Lev_{n+1}(T')$. Consider the function $G:T'\rightarrow X$ defined by $$G(s)=\begin{cases}
        0, & s\in T'\restriction \omega\setminus\{n\}\\
        \gamma_s, & s\in \Lev_n(T') \end{cases}.$$
 By the induction hypothesis there is $T^{(n+1)}\subseteq T'\subseteq T^{(n)}$ and $\gamma^*$ such that $$G\restriction \Lev_n(T^{(n+1)})\equiv\gamma^*.$$ We claim that $F\restriction \Lev_{n+1}(T^{(n+1)})\equiv\gamma^*$. To see this, fix $t^{\frown}\alpha\in \Lev_{n+1}(T^{(n+1)})$, then $\alpha\in A_t$ and $t\in \Lev_n(T^{(n+1)})$. Thus
$$F(t^{\frown}\alpha)=F_t(\alpha)=\gamma_\alpha=G(t)=\gamma^*.$$
This concludes the induction. Let $T^*=\bigcap_{n<\omega}T^{(n)}$, note that each level is an intersection of finitely many sets. Therefore, $T^*$ is a  $\vec{U}$-fat tree and $\forall n<\omega, \  F\restriction \Lev_n(T^*)$ is constant. $\square$

 Next we prove the Prikry condition.
\begin{theorem}
Assume that every filter in $\vec{U}$ is an ultrafilter. Let $\sigma$ be a formula in the forcing language and $p\in P_T(\vec{U})$. Then there exists $p\leq^* p^*$ such that $p^* \ || \sigma$.
\end{theorem}
\textit{Proof}:
Let $p\in P_T(\vec{U})$. Find $p\leq^*p^*$ so that for every $t\in T(p^*)$, $$p^{*\frown}\langle t\rangle || \sigma \  \Leftrightarrow \ \exists q\geq^* p^{*\frown}\langle t\rangle  \ q|| \sigma.$$ 
To find such $p^*$, let $t\in T(p)$. If $\exists q\geq^* p^{*\frown}\langle t\rangle$ such that  $ q|| \sigma$, set $S(q)=T(q)$, otherwise, set $S(t)=(T)_t$. By proposition \ref{PropTree}, amalgamate $\langle S(t)\mid t\in T(p)\rangle$ to a single tree $T^*\subset T(p)$ with the same trunk $t(p)$. Then the condition $p^*=\langle t(p),T^*\rangle$ is as wanted. Next,
 define a function $F:T(p)\rightarrow \{0,1,2\}$ by 
$$ F(t)=\begin{cases}
        0, &  p^{\frown}\langle t\rangle\Vdash \sigma\\
        1, &  p^{\frown}\langle t\rangle\Vdash \neg\sigma\\
        2, & \text{Otherwise } \end{cases}.$$
By theorem \ref{RowTree}, there is $T'$ homogeneous for $F$. Find minimal $n<\omega$ such that $F\restriction \Lev_n(T')$ is constantly $0$ or $1$, suppose for example it is constantly $0$. 
Toward a contradiction, suppose that $n>0$, pick any $t\in \Lev_{n-1}(T')$.
We claim that $F(t)=0$, this will contradict the minimality of $n$.
Otherwise, there is $p^{\frown}t\leq q$ such that $q\Vdash \neg \sigma$. If necessary, extend $q$ so that $p^{\frown}\langle t^{\frown}\alpha\rangle\leq q$, for some $\alpha\in suc_{T'}(t)$.
We have reached a contradiction since $$p^{\frown}\langle t^{\frown}\alpha\rangle\Vdash\sigma.$$
Thus $n=0$, which meas that $p^{\frown} \langle T'\rangle \ ||\ \sigma$.
$\square$
\begin{corollary}
Let $\nu<\kappa$ and assume that $\vec{U}$ is a tree of filters such that for every $a\in[\kappa]^{<\omega}$, $U_a$ is $\nu^+$-complete ultrafilter. Then $P^{V[G]}(\nu)=P^{V}(\nu)$.
\end{corollary}
\begin{corollary}
If for every $a\in[\kappa]^{<\omega}$, $U_a$ is a $\kappa$-complete ultrafilter, then $P_T(\vec{U})$ preserves cardinals.
\end{corollary}

Note that in the last two theorems we do not need every filter to be a $\kappa$-complete ultrafilter, but that for every $\vec{U}$-fat tree $T$ there is $T'\subseteq T$ such that for every $a\in T'$, $U_a$ is a $\kappa$-complete ultrafilter. We will prove that this requirement is necessary.
\begin{proposition} Let $\vec{U}$ be a tree of filters.
\begin{enumerate} 
\item For every $\lambda<\kappa$, if there exists $a^*\in[\kappa]^{<\omega}$ such that the set $$A:=\{a\in[\kappa]^{<\omega}\mid  U_a \text{\ is \ }\lambda-\mathrm{complete\ and\ not \ }\lambda^+-\mathrm{complete}\}$$ is dense above $a^*$ i.e. for every $\vec{U}$-fat tree $T$ with trunk above $a^*$, $A\cap T\neq\emptyset$. Then there is a generic extension of $P_T(\vec{U})$ which adds a new $\omega$-sequence to $\lambda$.
\item If there exists $a^*\in[\kappa]^{<\omega}$ such that the set $$B:=\{a\in[\kappa]^{<\omega}\mid U_a \mathrm{\  is \ not\ an\ ultrafilter}\}$$ is dense above $a^*$, there is a generic extension of $P_T(\vec{U})$ which adds a new real.
\end{enumerate}
\end{proposition}
\textit{Proof}:
For $(1)$, fix $p$ such that $t(p)=a^*$ and let $G$ be generic with $p\in G$. For every $a\in A$, fix $\langle X_{i,a}\mid i<\lambda\rangle\subseteq U_a$ such that $$\forall i<j<\lambda, X_{i,a}\supseteq X_{j,a} \text{ and } \cap_{i<\lambda}X_{i,a}\notin U_a.$$
Define in $V[G]$, $g:\omega\rightarrow \lambda$ by $$g(n)=\min(i<\lambda\mid C_G(n)\notin X_{i,C_G\restriction n})$$ and if $C_G(n)\in X_{i,C_G\restriction n}$ for every $i<\lambda$, define $g(n)=0$.
Claim that $g\notin V$, otherwise, let $p\leq q\in G$ be such that $q\Vdash \dot{g}=g$. By density of $A$, there is $ t\in A\cap T(q)$. 
Let $$g(|t|+1)=\gamma^*<\lambda,$$ then $X_{\gamma^*,t}\in U_t$. Find the minimal $j>\gamma^*$ such that $\succ_{T(q)}(t)\cap X_{t,\gamma^*}\setminus X_{t,j}\neq\emptyset$. 
There is such $j$ since $$\succ_{T(q)}(t)\cap X_{t,\gamma^*}\in U_t\text{ and }\cap_{i<\lambda}X_{i,t}\notin U_a.$$
Pick $x\in \succ_{T(q)}(t)\cap X_{t,\gamma^*}$ and extend $q$ to $q^*=q^{\frown}\langle t,x\rangle$. Then $$q^*\Vdash \dot{g}(|t|+1)=g(|t|+1)=\gamma^* \text{ and } q^*\Vdash \gamma^*< \dot{g}(|t|+1)=j,$$ contradiction. 
For $(2)$, let $a\in B$ and $X_a$ such that $X_a,\kappa\setminus X_a\notin U_a$.  
Define the real $$r=\{n<\omega\mid C_G(n)\in X_{C_G\restriction n}\}.$$
The proof that $r\notin V$ is completely analogous to the one given in \ref{newreal} for Prikry forcing.
$\square$
\begin{definition}
Let $G$ be $P_T(\vec{U})$-generic. Define  $C_G=\bigcup\{t(p)\mid p\in G\}$. A sequence $C\in[\kappa]^\omega$ is \textit{Tree-Prikry-generic for $\vec{U}$}, if $C=C_G$ for some generic $G\subseteq P_T(\vec{U})$. 
\end{definition}
\begin{lemma}
Let $\vec{U}$ be a tree of filters such that for every $a\in[\kappa]^{<\omega}$, $U_a$ contains all the final segments. Suppose $G$ is a generic filter for $P_T(\vec{U})$, then:
\begin{enumerate}
\item $\otp(C_G)=\omega$.
\item $\sup(C_G)=\kappa$.
\item For every $\langle A_a\mid a\in[\kappa]^{<
\omega}\rangle$ such that $A_a\in U_a$, there exists $n<\omega$ such that for every $m\geq n$, $ C_G(m)\in A_{C_G\restriction m}$.
\item For every $\langle T(a)\mid a\in[\kappa]^{<
\omega}\rangle$ such that $T(a)$ is a tree with trunk $a$, there exists $n<\omega$ such that for every $m\geq n$, $ C_G\in mb(T(C_G\restriction m))$.

\end{enumerate}
\end{lemma}
\textit{Proof:} For $(1)$, it suffices to prove that $C_G$ is infinite and $\forall \nu<\kappa$, $\nu\cap C_G$ is finite. This follows by the density of the sets:
$$D_n=\{p\in P_T(\vec{U})\mid n(p)\geq n\},  \ \text{ for } n<\omega.$$
$$E_\delta=\{p\in P_T(\vec{U})\mid \min(\succ_{T(p)}(t(p)))>\delta\}, \ \text{ for } \delta<\kappa.$$
$(2)$ follows by the density of the sets:
$$F_\delta=\{p\in P_T(\vec{U})\mid \max(t(p))>\delta\}, \ \text{ for } \delta<\kappa.$$
For $(3)$, let $p\in P_T(\vec{U})$.  By proposition \ref{PropTree} $(1)$, we can shrink $T(p)$
to $T^*$ so that for every $s\in T^*$, $\succ_{T^*}(s)\subseteq A_s$. Let $p\leq^*p^*$ defined by $$p^*=p^{\frown}\langle T^*\rangle.$$
By density, find such $p^*\in G$. It follows that for every $m> n(p^*)$,   $C_G(m)\in A_{C_G\restriction m}$.
Finally for $(4)$, let $p\in P_T(\vec{U})$. By proposition \ref{PropTree} $(3)$, it is possible to amalgamate the trees $T(a)$ to a tree $T^*$, so that $$\forall s\in T^*, \ (T^*)_s\subseteq T(s).$$ Let $p^*=p^{\frown}\langle T^*\rangle\geq^*p$. By density, find such $p^*\in G$. Let $m> n(p^*)$  and let $k<\omega$. There is $r\in G$ such that $|t(r)|>k$ and $p^{*\frown}\langle C_G(n+1),...,C_G(m-1)\rangle\leq r$, in particular, $$C_G\restriction |t(r)|=t(r)\in T(C_G\restriction m).$$ Since $T(C_G\restriction m)$ is downward closed, $C_G\restriction k\in T(p^*)$. Hence $\forall k<\omega$, $C_G\restriction k\in T(C_G\restriction m)$ and by definition of the set of maximal branches, $$C_G\in mb(T(C_G\restriction m)).$$
$\square$
\begin{corollary}
$cf^{V[G]}(\kappa)=\omega$.
\end{corollary}
\begin{definition}\label{FilterTree}
Let $C\in[\kappa]^{\omega}$, define the \textit{filter generated by $C$},
$$G_C=\{p\in P_T(\vec{U})\mid p(t)=C\cap (\max(p(t))+1)\wedge C\in mb(T(p))\}.$$
\end{definition}
\begin{lemma}
For every $C$, $G_C$ is a filter. If $C=C_G$ then $G_C=G$ and $G^*_C=G$.
\end{lemma}
\textit{Proof}: Analogous to lemma \ref{Filter}.
$\square$

Next we prove the strong Prikry property for Tree-Prikry forcing.
\begin{lemma}\label{strongPrikryTree}
Let $\vec{U}$ be a tree of filters such that for every $a\in[\kappa]^{<\omega}$, $U_a$ is an ultrafilter which contains all the final segments. Then for every $p\in P_T(\vec{U})$, there is $p^*\geq^*p$ and a natural number $n<\omega$ such that for every $t\in \Lev_n(T(p^*))$, $p^{*\frown}t\in D$.
\end{lemma}
\textit{Proof}:
Define $F:T(p)\rightarrow \{0,1\}$ by $$F(t)=0 \ \Leftrightarrow \ \exists  T'\subseteq (T(p))_t \ p^{\frown}\langle t,T'\rangle\in D.$$ 
By \ref{RowTree}, there is a $\vec{U}$-fat tree $T'\subseteq T(p)$ with trunk $t(p)$, such that  $F\restriction \Lev_n(T')$ is constant for every $n<\omega$. Let $t\in T'$, if $F(t)=0$ set $S(t)$ to be some witnessing tree for it, otherwise let $S(t)=(T')_t$. By \ref{PropTree} $(3)$,  find $T''$ such that 
$$\forall t\in T'' \ (T'')_t\subseteq S(t).$$
Define $T^*=T'\cap T''$. Since $D$ is dense, there is $s\in T^*$ such that for some $S$, $p^{\frown}\langle s,S\rangle\in D,$ hence $F(s)=0$. Let us claim that $p^*=p^{\frown}\langle T^*\rangle$ and $ht(s)$ are as wanted. Take any $s'\in \Lev_{ht(s)}(T^*)$, by homogeneity of $T'$, $F(s)=F(s')=0$. Therefore, $ p^{\frown}\langle s',T(s')\rangle\in D$. The tree $T^*$ has the property that $(T^*)_{s'}\subseteq T(s')$, thus $p^{\frown}\langle s', T(s')\rangle\leq p^{*\frown}\langle s'\rangle$. Finally, $p^{*\frown}\langle s'\rangle\in D$ since $D$ is open.
$\square$

The following is a Mathias-like characterization for Tree-Prikry forcing. 
\begin{theorem}\label{MathiasTree}
 Let $\vec{U}$ be a tree of filters such that for every $a\in[\kappa]^{<\omega}$, $U_a$ is an ultrafilter which contains all the final segments. If $C\in[\kappa]^{\omega}$ is a sequence such that:
 \begin{enumerate}
 \item $\sup(C)=\kappa$.
 \item For every $\langle T(a)\mid a\in[\kappa]^{<\omega}\rangle$ such that $T(a)$ is a $\vec{U}$-fat tree with trunk $a$, there exists $n<\omega$ such that for every $n\leq m<\omega$, $C\in mb(T(C\restriction m))$.
 \end{enumerate}
 Then $C$ is Tree-Prikry-generic for $\vec{U}$.
\end{theorem}
\textit{Proof}: Let $G=G_C$ be the filter generated by $C$ as defined in \ref{FilterTree}. Let $D$ be dense and open,
for every $s\in[\kappa]^{<\omega}$, apply the strong Prikry property to $D$ and the condition $\langle s, R(s)\rangle$, where $R(s)$ is the tree with trunk $s$ and for every $s\leq t\in R(s)$, $$\succ_{R(s)}(t)=\kappa\setminus(\max(t)+1).$$
Find a $\vec{U}$-fat tree $T(s)\subseteq R(s)$ and $n_s<\omega$ such that $$\forall t\in \Lev_{n_s}(T(s)) \ \langle t, (T(s))_t\rangle\in D.$$
We define $\langle S(a)\mid a\in[\kappa]^{<\omega}\rangle$ by induction on $|a|$. Let  $S(\langle\rangle)=T(\langle\rangle)$. Assume that $S(a)$ is defined and let  $\alpha\in\kappa\setminus\max(a)+1$. If $\alpha\in \succ_{S(a)}(a)$ define   $$S(a^{\frown}\langle\alpha\rangle)=(S(a))_{a^{\frown}\langle\alpha\rangle}\cap T(a^{\frown}\langle\alpha\rangle),$$ otherwise,  set $S({a^{\frown}\langle\alpha\rangle})=T({a^{\frown}\langle\alpha\rangle})$. 
By requirement (2) of the sequence $C$, $$\exists N<\omega \ \forall n\geq N \ C\in mb(S( C\restriction n)).$$ Consider the condition
$p=\langle C\restriction N,S(C\restriction N)\rangle\in P_T(\vec{U})$ and let $n^*=n_{C\restriction N}$. Then $$C\restriction n^*\in \Lev_{n^*}(S(C\restriction N))\subseteq \Lev_{n^*}(T(C\restriction N)).$$
We claim that $p^*=\langle C\restriction n^*,S(C\restriction n^*)\rangle\in D\cap G$. Since $n^*\geq N$, $C\in mb(S(C\restriction n^*))$ and by definition of $G$, $p^*\in G$. To see that $p^*\in D$, note that for $n\geq N$, 
$C\restriction n\in S(C\restriction N)$ and by the recursive definition of $S$, $$S(C\restriction n)=(S(C\restriction N))_{C\restriction n}\subseteq (T(C\restriction N))_{C\restriction n}.$$
Particularly, for $n^*$ 
$$p^*=\langle C\restriction n^*,S(C\restriction n^*)\rangle=\langle C\restriction n^*,(S(C\restriction N))_{C\restriction n^*}\rangle\geq \langle C\restriction n^*, 
(T(C\restriction N))_{C\restriction n^*}\rangle\in D.$$
Since $D$ is open, $p^*\in D$.
$\square$
\begin{corollary}\label{ImpMathiasTree}
 Let $\vec{U}$ be a tree of filters such that for every $a\in[\kappa]^{<\omega}$, $U_a$ is an ultrafilter which contains all the final segments. If $C\in[\kappa]^{\omega}$ is a sequence such that:
\begin{enumerate}

\item $\sup(C)=\kappa$.
 \item For every $\langle A_a\mid a\in[\kappa]^{<\omega}\rangle$ such that $A_a\in U_a$, there exists $n<\omega$ such that for every $n\leq m<\omega$,
 $C(m)\in A_{C\restriction m}$.
 \end{enumerate}
 Then $C$ is Tree-Prikry-generic for $\vec{U}$.
\end{corollary}
\textit{Proof}: Prove condition (2) of theorem \ref{MathiasTree}. Let $\langle T(a)\mid a\in [\kappa]^{<\omega}\rangle$ be a tree of trees and define sets $A_a$ for $a\in[\kappa]^{<\omega}$. $A_{\langle\rangle}=\succ_T(\langle\rangle)$, let $i\leq |a|$ be minimal such that for every $i\leq j\leq |a|$, $a\in T(a\restriction j)$ and set $$A_a=\bigcap_{i\leq j\leq |a|}\succ_{T(a\restriction j)}(a)\in U_a.$$ By propery $(2)$ of $C$, there is $N$ such that for every $n\geq N$, $C(n)\in A_{C\restriction n}$. Fix $N\leq k<\omega$, let us prove that $C\in mb(T(C\restriction k))$. By induction on $m\geq k$, prove that $C\restriction m\in T(C\restriction k)$. For $m=k$ this is just by the fact that $C\restriction k$ is the trunk of $T(C\restriction k)$. Assume that $C\restriction m\in T(C\restriction k)$, then $$\succ_{T(C\restriction k)}(C\restriction m)\supseteq A_{C\restriction m}.$$ 
Moreover, $C(m)\in A_{C\restriction m}$ implying that $C(m)\in \succ_{T(C\restriction k)}(C\restriction m)$. Finally, we conclude that $C\restriction (m+1)\in T(C\restriction k)$ and the induction step follows.
$\square$

Next we translate this characterization to Tree-Prikry with a $\kappa$-sequence, Tree-Prikry with an $ \omega$-sequence and Tree-Prikry with a single ultrafilter. This criteria was asked for in \cite{MinimalPrikry}.

\begin{corollary} Let $\vec{W}=\langle W_\alpha\mid \alpha<\kappa\rangle$ be a $\kappa$-sequence of $\kappa$-complete ultrafilters. Let $C\in[\kappa]^{\omega}$, then $C$ is Tree-Prikry-generic for $\vec{W}$ iff: 
\begin{enumerate}
    \item $\sup(C)=\kappa$.
    \item For every $\langle A_\alpha\mid a<\kappa\rangle$ such that $A_\alpha\in W_\alpha$,  there exists $n$ such that for every $m\geq n$, $C(m)\in A_{C(m-1)}$.
\end{enumerate}
\end{corollary}
\textit{Proof}:
     Suppose $C$ satisfies $$(*) \ \ \forall\langle A_\alpha\mid a<\kappa\rangle \ \exists n \  \forall m\geq n \ C(m)\in A_{C(m-1)}.$$ Let us show the criteria in \ref{ImpMathiasTree} holds. Fix $\langle A_a\mid a\in[\kappa]^{<\omega}\rangle$ such that $A_a\in U_a$, define $$A_\alpha=\bigcap_{a\in[\alpha]^{<\omega}}A_{a^{\frown}\langle\alpha\rangle}.$$ Since $W_\alpha$ is $\kappa$-complete, $A_\alpha\in W_\alpha$. By $(*)$, there is $m$ such that for every $n>m$, $$C(n)\in A_{C(n-1)}\subseteq A_{C\restriction n}.$$
     $\square$
     \begin{corollary}\label{MathiasOmega} Let $\vec{W}=\langle W_n\mid n<\omega\rangle$ be an $\omega$-sequence of $\kappa$-complete ultrafilters and let $C\in[\kappa]^{\omega}$. For every $n<\omega$, fix $\pi_n:\kappa\rightarrow \kappa$ such that $[\pi_n]_{W_n}=\kappa$. Then $C$ is Tree-Prikry-generic for $\vec{W}$ iff:
     \begin{enumerate}
     \item $\sup(C)=\kappa$.
         \item For every $\langle A_n\mid n<\omega\rangle$ such that $A_n\in U_n$, there exists $n$ such that for $m\geq n$, $C(m)\in A_{m-1}$.
         \item There exists $n$ such that for every $m\geq n$, $\pi_{m-1}(C(m))>C(m-1)$.
         \end{enumerate}
     \end{corollary}
    \textit{Proof}: Suppose $C_G$ is generic for $P_T(\vec{W})$. Let $p\in P_T(\vec{W})$, for every $t\in T(p)$ define
    $$X_t=\{\nu<\kappa\mid \pi_{|t|}(\nu)>\max(t)\}.$$
     To see that $X_t\in W_{|t|}$, note that in the ultrapower by $W_{|t|}$, $$\max(t)=j_{W_{|t|}}(\max(t))<\kappa=[\pi_{|t|}]_{W_{|t|}}.$$ Shrink $\succ_{T(p)}(t)$ to $\succ_{T(p)}(t)\cap X_t$. By density, find $p^*\in G$ such that for every $t\in T(p^*)$, $\succ_{T(p^*)}(t)\subseteq X_t$. Then $C_G\in mb(T(p^*))$ and $$\forall m> n(p^*) \ \pi_{m-1}(C_G(m))>C_G(m-1).$$ For the other direction, we prove the condition in \ref{ImpMathiasTree}, let $\langle A_a\mid a\in[\kappa]^{<\omega}\rangle$ such that $A_a\in W_{|a|}$. Define $A_0=A_{\langle\rangle}\in W_0$ and 
    $$A_{n}=\Delta^*_{a\in[\kappa]^n}A_a:=\{\nu<\kappa\mid \forall a\in[\pi_{n}(\nu)]^{n} \  \nu\in A_a\}.$$
    To see that $A_{n}\in W_{n}$, note that in the ultrapower by $W_n$, $\kappa=[\pi_n]_{W_n}$ and $$\forall a\in[[\pi_n]_{W_n}]^{n}, \ [id]_{W_n}\in j_{W_n}(A_a).$$ By L{\'o}s theorem $A_n\in W_n$. There is $m<\omega$ such that for every $n>m$, $$\pi_{n-1}(C(n))>C(n-1)\text{ and }C(n)\in A_{n-1}.$$ Since $C\restriction n\in [\pi_{n-1}(C(n))]^n$, $C(n)\in A_{C\restriction n}$.
$\square$
\begin{corollary}\label{MathiasSingle} Let $U$ be a $\kappa$-complete ultrafilter and let $C\in[\kappa]^{\omega}$. Fix $\pi:\kappa\rightarrow\kappa$ such that $[\pi]_U=\kappa$. Then $C$ is Tree-Prikry-generic for $U$ iff:
\begin{enumerate}
\item $\sup(C)=\kappa$.
\item  For every $A\in U$, there exists $n$ such that for every $m\geq n$, $C(m)\in A$.
\item There exists $n$ such that for every $m\geq n$, $\pi(C(m))>C(m-1)$.
\end{enumerate}
\end{corollary}
\textit{Proof}: To prove condition $(2)$ of corollary \ref{MathiasOmega}, given $\langle A_n\mid n<\omega\rangle$ such that $A_n\in U$, apply condition $(2)$ to the set $A=\cap_{a<\omega}A_n\in U$.
$\square$

In corollaries \ref{MathiasOmega} and \ref{MathiasSingle}, we required that the sequence satisfies  $\pi_{n-1}(C(n))>C(n-1)$. This condition emphasis the difference between Prikry-generic sequences and Tree-Prikry-generic sequences, it was crucial to diagonally intersect sets of the same level. More precisely, the notion of Tree-Prikry generic sequences, generalizes of regular Prikry-generic sequences since if $U$ is normal, the new requirement is satisfied automatically as $[id]_U=\kappa$.

We wish to generalize the Solovay observation that the critical sequence of the iterated ultrapower by a normal measure $U$ is Prikry-generic for $P(j_{\omega}(U))$  above $M_\omega$. Let $\vec{U}$ be a tree of $\kappa$-complete ultrafilters over $\kappa$, denote $V=M_0$. Fix the first measure in $\vec{U}$ and its embedding i.e. $$U_0:=U_{\langle\rangle}, \ \kappa_0:=\kappa, \ \vec{U}_0:=\vec{U}, \ 
j_0:V\rightarrow Ult(V,U_0).$$
Set $M_1:=Ult(V,U_0)$, $\kappa_1:=j_0(\kappa), \ \vec{U}_1:=j_0(\vec{U})$. Consider $\delta_0:=[id]_{U_0}$, take the ultapower by $$U_1:=(\vec{U}_1)_{\langle \delta_0\rangle}, \ 
j_1:M_1\rightarrow Ult(M_1,U_1).$$
Once again, denote $M_2:=Ult(M_1,U_1), \ \kappa_2:=j_1(\kappa_1),\ \vec{U}_2:=j_1(\vec{U}_1)$.
Generally, assume we have defined $$M_n,\ \vec{U}_n,\ \kappa_n, \ U_n=(\vec{U}_n)_{\langle\delta_0,...\delta_{n-1}\rangle}, \ j_{n}:M_n\rightarrow Ult(M_n,U_n).$$
Set $$\delta_n:=[id]_{U_n}, \  M_{n+1}:=Ult(M_n,U_n),\ \kappa_{n+1}:=j_{n}(\kappa_n), \  \vec{U}_{n+1}:=j_{n}(\vec{U}_n),\ U_{n+1}=(\vec{U}_{n+1})_{\langle\delta_0,...,\delta_n\rangle}.$$ For $n<m$, let $$j_{n,m}=j_{m-1}\circ...\circ j_n:M_n\rightarrow M_m,$$ and if $n=m$ let $j_{n,n}=id$. The system of models $\langle M_n,j_{n,m}\mid n<m<\omega\rangle$ has a direct limit denoted by $M_\omega$ and direct limit embeddings $$j_{n,\omega}:M_n\rightarrow M_{\omega}, \ j_{n,\omega}=j_{m,\omega}\circ j_{n,m},$$ such that $M_\omega=\bigcup_{n<\omega}j_{n,\omega}[M_n]$. The model $M_{\omega}$ is well founded and we identify it with its transitive collapse.
\begin{theorem}
Let $\vec{U}$ be a tree of $\kappa$-complete ultrafilters. The sequence $\langle\delta_n\mid n<\omega\rangle$ is Tree-Prikry-generic over $M_\omega$ for $\vec{U}_{\omega}:=j_{0,\omega}(\vec{U})$.
\end{theorem}
\textit{Proof}: We verify that the Mathias criteria holds for the sequence  $\langle\delta_n\mid n<\omega\rangle$. Let $$\vec{A}=\langle A_a\mid a\in [\kappa_{\omega}]^{<\omega}\rangle\in M_{
\omega}, \ A_a\in (\vec{U}_{\omega})_a.$$ Then there is $n<\omega$ and $x\in M_n$ such that $j_{n,\omega}(x)=\vec{A}$. By elementarity, $x$ is also a tree of filters, $x=\langle B_a\mid a\in[\kappa_n]^{<\omega}\rangle$.
We claim that for $m\geq n$, $\delta_m\in A_{\langle \delta_0,...,\delta_{m-1}\rangle}$. Let $Y=j_{n,m}(x)$, then $$Y_{\langle\delta_0,...,\delta_{m-1}\rangle}\in (\vec{U}_m)_{\langle\delta_0,...,\delta_{m-1}\rangle}=U_m.$$
It follows that $\delta_m=[id]_{U_m}\in j_{m,m+1}(Y_{\delta_0,...,\delta_{m-1}})$, thus by elementarity and the fact that for $i<m$, $\delta_i<\kappa_m=crit(j_{m,\omega})$
$$\delta_m=j_{m+1,\omega}(\delta_m)\in j_{m+1,\omega}(j_{m,m+1}(Y_{\delta_0,...,\delta_{m-1}})).$$
Hence
$$\delta_m \in j_{m,\omega}(Y)_{\delta_0,...,\delta_{m-1}}=j_{n,\omega}(x)_{\delta_0,...,\delta_{m-1}}=A_{\delta_0,...,\delta_{m-1}}.$$
$\square$

One major difference between Tree-Prikry forcing and the Prikry forcing is the structure of intermediate models. In the case of Prikry forcing with we always have non-trivial intermidiate extensions, it is not always the case for Tree-Prikry forcing. In \cite{MinimalPrikry}, it is proven that in the case of a $\kappa$-sequence of distinct normal measures, or any $\kappa$-sequence of measures that can be separated, we do not have non-trivial intermediate extensions.
\begin{theorem}
Let $\vec{U}=\langle U_i\mid i<\kappa\rangle$ be a $\kappa$-sequence of distinct normal ultrafilters. Let  $G\subseteq P_T(\vec{U})$ be generic over $V$, then every $A\in V[G]$ either $A\in V$ or $V[A]=V[G]$.
\end{theorem}
\section{Questions}
To the best of our knowledge the following questions are open:
\begin{question}
Is there a combinatorical characterization of all $\vec{U}$ tree of filters for which $P_T(\vec{U})$ preserve cardinals?
\end{question}
\begin{question}
Is there a combinatorical characterization of all $\vec{U}$ tree of filters for which $P_T(\vec{U})$ has the mathias characterization? does the answer changes according to $\kappa$-sequence/$\omega$-sequence/single ultrafilter?
\end{question}
\begin{question}
What are all intermediate models of $P_T(\vec{U})$?
\end{question}
\begin{question}
Let $U$ be a non-normal ultrafilter, how does the intermediate models of $P(U)$ looks like?
\end{question}

\end{document}